\magnification=1200

\loadmsam
\loadmsbm
\loadeufm
\loadeusm
\UseAMSsymbols

\font\BIGtitle=cmr10 scaled\magstep3
\font\bigtitle=cmr10 scaled\magstep1
\font\boldsectionfont=cmb10 scaled\magstep1
\font\section=cmsy10 scaled\magstep1

\def\scr#1{{\fam\eusmfam\relax#1}}

\def\scrA{{\scr A}}
\def\scrB{{\scr B}}
\def\scrC{{\scr C}}
\def\scrD{{\scr D}}
\def\scrE{{\scr E}}
\def\scrF{{\scr F}}

\def\scrL{{\scr L}}
\def\scrK{{\scr K}}

\def\scrM{{\scr M}}
\def\scrN{{\scr N}}

\def\scrT{{\scr T}}

\def\scrW{{\scr W}}
\def\gr#1{{\fam\eufmfam\relax#1}}

\def\grD{{\gr D}}

	\def\grg{{\gr g}}

\def\grL{{\gr L}}

\def\grO{{\gr O}}	
	\def\grp{{\gr p}}

\def\grS{{\gr S}}	
\def\grT{{\gr T}}	
	
\def\grV{{\gr V}}

\def\db#1{{\fam\msbfam\relax#1}}

\def\dbA{{\db A}} 
\def\dbC{{\db C}} 
 \def\dbF{{\db F}}
\def\dbG{{\db G}} \def\dbH{{\db H}}

 \def\dbN{{\db N}}
 
\def\dbQ{{\db Q}} \def\dbR{{\db R}}
\def\dbS{{\db S}}

 \def\dbZ{{\db Z}}

\def\der{\text{der}}
\def\Sh{\hbox{\rm Sh}}

\def\sp{\text{sp}}
\def\Res{\text{Res}}

\def\ad{\text{ad}}
\def\Ad{\text{Ad}}

\def\Gal{\text{Gal}}
\def\Hom{\text{Hom}}
\def\End{\text{End}}
\def\Spec{\text{Spec}}
\def\gen{\text{gen}}

\def\sc{\text{sc}}

\def\Lie{\text{Lie}}

\def\leaderfill{\leaders\hbox to 1em
     {\hss.\hss}\hfill}
\def\nspace{\lineskip=1pt\baselineskip=12pt\lineskiplimit=0pt}

\def\finishproclaim{\par\rm
     \ifdim\lastskip<\medskipamount\removelastskip
     \penalty55\medskip\fi}
\def\endproof{$\hfill \square$}
\def\proof{\par\noindent {\it Proof:}\enspace}
\def\references#1{\par
  \centerline{\boldsectionfont References}\medskip
     \parindent=#1pt\nspace}
\def\Ref[#1]{\par\hang\indent\llap{\hbox to\parindent
     {[#1]\hfil\enspace}}\ignorespaces}
\def\Item#1{\par\smallskip\hang\indent\llap{\hbox to\parindent
     {#1\hfill$\,\,$}}\ignorespaces}
\def\ItemItem#1{\par\indent\hangindent2\parindent
     \hbox to \parindent{#1\hfill\enspace}\ignorespaces}

\def\Le{{\mathchoice{\,{\scriptstyle\le}\,}
  {\,{\scriptstyle\le}\,}
  {\,{\scriptscriptstyle\le}\,}{\,{\scriptscriptstyle\le}\,}}}
\def\Ge{{\mathchoice{\,{\scriptstyle\ge}\,}
  {\,{\scriptstyle\ge}\,}
  {\,{\scriptscriptstyle\ge}\,}{\,{\scriptscriptstyle\ge}\,}}}

\def\arrowsim{\,\smash{\mathop{\to}\limits^{\lower1.5pt
  \hbox{$\scriptstyle\sim$}}}\,}

\def\doublemaprights#1#2#3#4{\raise3pt\hbox{$\mathop{\,\,\hbox to     
#1pt{\rightarrowfill}\kern-30pt\lower3.95pt\hbox to
     #2pt{\rightarrowfill}\,\,}\limits_{#3}^{#4}$}}

\def\rightcapdownarrow{\raise9pt\hbox{$\ssize\cap$}\kern-7.75pt
     \Big\downarrow}

\def\rcapmapdown#1{\rightcapdownarrow\kern-1.0pt\vcenter{
     \hbox{$\scriptstyle#1$}}}

\def\rmapdown#1{\Big\downarrow\kern-1.0pt\vcenter{
     \hbox{$\scriptstyle#1$}}}
\def\rightsubsetarrow#1{{\ssize\subset}\kern-4.5pt\lower2.85pt
     \hbox to #1pt{\rightarrowfill}}
\def\longtwoheadedrightarrow#1{\raise2.2pt\hbox to #1pt{\hrulefill}
     \!\!\!\twoheadrightarrow}

\def\Gal{\operatorname{\hbox{Gal}}}
\def\Hom{\operatorname{\hbox{Hom}}}

\def\im{\hbox{Im}}

\NoBlackBoxes
\parindent=25pt
\document
\footline={\hfil}

\noindent 
\centerline{\BIGtitle Projective integral models of Shimura varieties}
\medskip
\centerline{\BIGtitle of Hodge type with compact factors}
\vskip 0.3 cm
\centerline{\bigtitle Adrian Vasiu, U of Arizona, final version 12/14/06, to appear in Crelle}
\footline={\hfill}
\vskip 0.3 cm
\noindent
{\bf ABSTRACT.} Let $(G,X)$ be a Shimura pair of Hodge type such that $G$ is the Mumford--Tate group of some elements of $X$. We assume that for each simple factor $G_0$ of $G^{\ad}$ there exists a simple factor of $G_{0\dbR}$ which is compact. Let $N\Ge 3$. We show that for many compact open subgroups $K$ of $G(\dbA_f)$, the Shimura variety $\Sh(G,X)/K$ has a projective integral model $\scrN$ over $\dbZ[{1\over N}]$ which is a finite scheme over a certain Mumford moduli scheme $\scrA_{g,1,N}$. Equivalently, we show that if $A$ is an abelian variety over a number field and if the Mumford--Tate group of $A_{\dbC}$ is $G$, then $A$ has potentially good reduction everywhere. The last result represents significant progress towards the proof of a conjecture of Morita. If $\scrN$ is smooth over $\dbZ[{1\over N}]$, then it is a N\'eron model of its generic fibre. In this way one gets in arbitrary mixed characteristic, the very first examples of general nature of projective N\'eron models whose generic fibres are not finite schemes over abelian varieties.   
\medskip\noindent
{\bf Key words}: abelian and semiabelian schemes, Mumford--Tate groups, Shimura varieties, integral models, and N\'eron models.
\medskip\noindent
{\bf MSC 2000}: Primary 11G10, 11G18, 14F30, 14G35, 14K10, 14K15, and 14J20.

\footline={\hss\tenrm \folio\hss}
\pageno=1

\bigskip\smallskip
\noindent
{\boldsectionfont \S1. Introduction}
\bigskip

Let $\dbS:=\Res_{\dbC/\dbR} \dbG_{m\dbC}$ be the two dimensional torus over $\dbR$ with the property that $\dbS(\dbR)$ is the  multiplicative group of non-zero complex numbers. Let $E$ be a number field. Let $O_E$ be the ring of integers of $E$. We fix an embedding $i_E:E\hookrightarrow\dbC$. Let $A$ be an abelian variety over $E$. Let $W_A:=H_1(A({\dbC}),\dbQ)$ be the first Betti homology group of the complex manifold $A(\dbC)$ with coefficients in $\dbQ$. Let $h_A\colon\dbS\to \text{GL}_{W_A\otimes_{\dbQ} \dbR}$
be the homomorphism that defines the Hodge $\dbQ$--structure on $W_A$. Let $H_A$ be the Mumford--Tate group of $A_{\dbC}$. We recall that $H_A$ is a reductive group over $\dbQ$ and that $H_A$ is the smallest subgroup of $\text{GL}_{W_A}$ with the property that $h_A$ factors through $H_{A\dbR}$, cf. [8], Propositions 3.6 and 3.4. Let $H_A^{\ad}$ be the adjoint group of $H_A$ i.e., the quotient of $H_A$ through its center.

\bigskip\noindent
{\bf 1.1. Definition.} We say the abelian variety $A$ has {\it compact factors}, if for each simple factor $H_0$ of $H_A^{\ad}$ there exists a simple factor of $H_{0\dbR}$ which is compact. 

\medskip
It is conjectured (see [20], page  437 and [26], Conjecture 3.1.3) that if the $\dbQ$--rank of $H_A^{\ad}$ is $0$ (i.e., if $\dbG_{m\dbQ}$ is not a subgroup of $H_A^{\ad}$), then there exists a finite field extension $E_1$ of $E$ such that $A_{E_1}$ extends to an abelian scheme over $O_{E_1}$ (i.e., such that $A_{E_1}$ has good reduction with respect to all finite primes of $E_1$). See [26] for other equivalent forms of this conjecture. Following [26], below we will refer to this conjecture as the Morita conjecture. We note down that if $A$ has compact factors, then each simple factor of $H_A^{\ad}$ has $\dbQ$--rank $0$ and thus $H_A^{\ad}$ itself has $\dbQ$--rank $0$. The goals of this paper are: (i) to prove the Morita conjecture under the assumption that $A$ has compact factors, and (ii) to reformulate and apply this result to integral models of Shimura varieties of Hodge type. 

\bigskip\noindent
{\bf 1.2. Basic Theorem.} {\it Suppose that $A$ has compact factors. Then there exists a finite field extension $E_1$ of $E$ such that $A_{E_1}$ extends to an abelian scheme over $O_{E_1}$.}

\medskip
Morita conjecture was first checked in some cases involving abelian varieties of PEL type (see [20]; see also [14], end of \S5). We recall that $A$ is of PEL type, provided $A_{\dbC}$  has a polarization $\lambda$ such that the derived group of $H_A$ is also the derived group of the intersection of $\text{GSp}(W_A,\psi)$ with the double centralizer of $H_A$ in $\text{GL}_{W_A}$ (here $\psi$ is the non-degenerate alternating form on $W_A$ defined by $\lambda$ and PEL stands for polarization, endomorphisms, and level structures). Some of the cases presented in [20] are not covered by the Basic Theorem. 

New cases of the validity of the Morita conjecture are provided in [24] and [26]. For instance, Paugam proved the Morita conjecture provided there exists a prime $p\in\dbN$ such that the $\dbQ_p$-rank of $H_{A\dbQ_p}^{\ad}$ is $0$ (see Lemma 2.3.1). The results [26], Propositions 4.2.2, 4.2.4, and 4.2.10 cover [24]; these results are also either covered by Lemma 2.3.1 or involve special cases when there exists a good prime $p\in\dbN$ for which a certain combinatorial condition on the natural action of $\Gal(\dbQ_p)$ on the set of simple factors of $H_{A\overline{\dbQ_p}}^{\ad}$ holds (see [26], Example 4.2.11 for such a concrete special case). Such good primes $p$ exist only if $A$ has compact factors and each simple factor $H_0$ of $H_A^{\ad}$ is ``simple enough" (like when $H_{0\dbR}$ has only one simple, non-compact factor). Good primes $p$ do not exist if there exists one simple factor $H_0$ of $H_A^{\ad}$ such that the following two properties hold: (i) $H_{0\dbC}$ is not isotypic of $A_n$ Lie type and (ii) the group $H_{0\dbR}$ either (ii.a) has more simple, non-compact factors than simple, compact factors or (ii.b) it is a Weil restriction $\Res_{F_0/\dbQ} \tilde H_0$, where $\tilde H_0$ is an absolutely simple, adjoint group over some ``arithmetically complicated" totally real number field $F_0$. Thus the results [26], Propositions 4.2.2, 4.2.4, 4.2.10, and 4.2.13 are particular cases of the combination of the Basic Theorem and Lemma 2.3.1.  

Basic Theorem, Lemma 2.3.1, and (the situations that can be reduced to) [14], end of \S5, form all (general) cases in which the Morita conjecture is presently known to hold. Lemma 2.3.1 and [14], end of \S5 pertain only to abelian varieties $A$ over $E$ for which, up to a replacement of $E$ by a finite field extension of it, there exists an abelian variety $B$ over $E$ that is of PEL type and that has the property that the three adjoint groups $H_A^{\ad}$, $H_{B}^{\ad}$, and $H_{A\times_E B}^{\ad}$ are isomorphic (cf. Remark 2.3.2 and [31], Corollary 4.10). 

The methods we use to prove the Basic Theorem pertain to Shimura varieties and rely on the constructions of [7], Proposition 2.3.10 and [29], Subsections 6.5 and 6.6. To outline the methods, in this paragraph we assume that $A$ has compact factors. We use the mentioned constructions in order to show that for a given prime $p\in\dbN$, up to a replacement of $E$ by a finite field extension of it, there exists an abelian variety $B$ over $E$ which has the following two properties: (a) the three adjoint groups $H_A^{\ad}$, $H_{B}^{\ad}$, and $H_{A\times_E B}^{\ad}$ are isomorphic, and (b) the monomorphism $H_B\hookrightarrow \text{GL}_{W_B}$ is ``manageable enough" so that we can check based on [10] and [9] that $B$ has good reduction with respect to all primes of $E$ that divide $p$ (see Section 3 and Subsection 4.1). Due to (a) and the good reduction part of (b), there exists a finite field extension $E_p$ of $E$ such that $A_{E_p}$ has good reduction with respect to all primes of $E_p$ that divide $p$ (see [26], Proposition 4.1.2; see also Subsection 4.1). As the field $E_1$ mentioned in the Basic Theorem, we can take the composite field of a suitable finite set of such fields $E_p$ (see Subsection 4.1).

Such ``extra" abelian varieties $B$, were first considered in [31] in connection to the Mumford--Tate conjecture for $A$ and in [26] in connection to the Morita conjecture for $A$. 

In Section 2 we review basic properties of Shimura varieties. In Section 3 we include a construction that is the very essence of the proof of the Basic Theorem and which in fact proves the Basic Theorem under some additional hypotheses. In Subsection 4.1 we prove the Basic Theorem. Example 4.2 is completely new. Basic Corollary 4.3 reformulates the Basic Theorem in terms of integral models of Shimura varieties of Hodge type. In Subsection 4.4 we apply the Basic Corollary 4.3 to provide new examples of general nature of N\'eron models in the sense of [3], page 12. In Example 4.5 and Remark 4.6 we apply the Basic Corollary 4.3 to correct an error in [29], Remark 6.4.1.1 2) and thus implicitly in [29], Subsubsection 6.4.11, for the cases to which the Basic Theorem applies (i.e., for Shimura pairs of preabelian type that have compact factors in a sense analogous to Definition 1.1; see Subsection 2.2 and Remark 4.6 (b) for precise definitions).

\smallskip\noindent
{\bf Acknowledgment.}  We would like to thank U of Arizona for providing us with good conditions with which to write this work. We would also like to thank the referee for many valuable suggestions and comments.

\bigskip\smallskip
\noindent
{\boldsectionfont \S2. Preliminaries}
\bigskip 

In Subsection 2.1 we gather different notations to be often used in the rest of the paper. See Subsection 2.2 for generalities on Shimura pairs and varieties. See Subsection 2.3 for simple properties that pertain to polarizations of $A$ and to the reductive group $H_A$. See Subsection 2.4 for the Shimura pairs naturally associated to $A$. See Subsection 2.5 for smooth toroidal compactifications. See Fact 2.6 and Proposition 2.7 for two results that pertain to an equivalent form of Theorem 1.2 to be stated explicitly in the Corollary 4.3. Lemma 2.8 will be used in Section 3.

\bigskip\noindent
{\bf 2.1. Notations and conventions.} A reductive group $H$ over a field $k$ is assumed to be connected. Let $Z(H)$, $H^{\der}$, and $H^{\ad}$ denote the center, the derived group, and the adjoint group (respectively) of $H$. We have $H^{\ad}=H/Z(H)$. Let $Z^0(H)$ be the maximal subtorus of $Z(H)$. Let $H^{\sc}$ be the simply connected semisimple group cover of $H^{\der}$. See [3], Subsection 7.6 for the Weil restriction of scalars functor $\Res_{k_1/k}$, where $k_1$ is a finite, \'etale $k$-algebra. We recall that if $H_1$ is a reductive group over $k_1$, then $\Res_{k_1/k} H_1$ is a reductive group over $k$ uniquely determined up to isomorphism by the group identities $\Res_{k_1/k} H_1(\heartsuit)=H_1(\heartsuit\otimes_k k_1)$ which are functorial on commutative $k$-algebras $\heartsuit$. For a free module $M$ of finite rank over a commutative ring with unit $R$, let $M^*:=\Hom_R(M,R)$ and let $\text{GL}_M$ be the group scheme over $R$ of linear automorphisms of $M$. If $\psi$ is a perfect alternating form on $M$, then $\text{GSp}(M,\psi)$ is viewed as a reductive group scheme over $R$. If $*$ or $*_R$  (resp. $*_+$ with $+$ as an arbitrary index different from $R$) is either an object or a morphism of the category of $\Spec(R)$-schemes, let $*_U$ (resp. $*_{+U}$) be its pull back via an affine morphism $\Spec(U)\to\Spec(R)$. Let $\overline{\dbQ}$ be the algebraic closure of $\dbQ$ in $\dbC$. We always use the notations of the first paragraph of Section 1. 

\bigskip\noindent
{\bf 2.2. On Shimura pairs.} A Shimura pair $(G,X)$ consists of a reductive group $G$ over $\dbQ$ and a $G(\dbR)$-conjugacy class $X$ of homomorphisms $\dbS\to G_{\dbR}$ that satisfy  Deligne's axioms of [7], Subsubsection 2.1.1: the Hodge $\dbQ$--structure on $\Lie(G)$ defined by any $x\in X$ is of type $\{(-1,1),(0,0),(1,-1)\}$, $\Ad\circ x(i)$ defines a Cartan involution of $\Lie(G^{\ad}_{\dbR})$, and no simple factor of $G^{\ad}$ becomes compact over $\dbR$. Here $\Ad:G_{\dbR}\to \text{GL}_{\Lie(G^{\ad}_{\dbR})}$ is the adjoint representation. These axioms imply that $X$ has a natural structure of a hermitian symmetric domain (see [7], Corollary 1.1.17). Similarly to Definition 1.1, we say $(G,X)$ has compact factors, if for each simple factor $G_0$ of $G^{\ad}$ there exists a simple factor of $G_{0\dbR}$ which is compact (thus $G_{0\dbR}$ has at least one simple, non-compact factor and at least one simple, compact factor). We note down that if $G$ is a torus (i.e., if the adjoint group $G^{\ad}$ is trivial), then $(G,X)$ has compact factors. For generalities on Shimura pairs and varieties and on their types, we refer to [6], [7], [17], [18], and [29], Subsections 2.2 to 2.5. To $(G,X)$ it is naturally associated a number field $E(G,X)$, called the reflex field of $(G,X)$ (see [6], [7], and [16]; see also the definition of the reflex field $E_A$ in the below Subsection 2.4). 

Let $\dbA_f:=\widehat\dbZ\otimes_{\dbZ}\dbQ$ be the ring of finite ad\`eles of $\dbQ$. For $K$ a compact open subgroup of $G(\dbA_f)$, let $\Sh(G,X)_{\dbC}/K:=G(\dbQ)\backslash (X\times G(\dbA_f)/K)$; it is a finite disjoint union of quotients of a (fixed) connected component of $X$ by arithmetic subgroups of $G(\dbQ)$. A theorem of Baily and Borel says that $\Sh(G,X)_{\dbC}/K$ has a canonical structure of a normal, quasi-projective ${\dbC}$-scheme (see [1], Theorem 10.11) which is smooth if $K$ is sufficiently small. Thus the projective limit $\Sh(G,X)_{\dbC}$ of the normal, quasi-projective $\dbC$-schemes $\Sh(G,X)_{\dbC}/K$'s, has a canonical structure of a regular $\dbC$-scheme. The Shimura variety $\Sh(G,X)$ is identified with the canonical model over $E(G,X)$ of the $\dbC$-scheme $\Sh(G,X)_{\dbC}$ (see [6], [7], [16], [18], and [19]). We have a natural right action of $G(\dbA_f)$ on the $E(G,X)$-scheme $\Sh(G,X)$ which is continuous in the sense of [7], Subsubsection 2.7.1. In particular, the quotient $\Sh(G,X)/K$ of $\Sh(G,X)$ through $K$ is a normal, quasi-projective $E(G,X)$-scheme which is ``the best arithmetic" model of $\Sh(G,X)_{\dbC}/K$.

\bigskip\noindent
{\bf 2.3. Simple properties.} The notion of good reduction of an abelian variety over the field of fractions of a discrete valuation ring, is stable under isogenies. Up to a replacement of $E$ by a finite field extension of it, $A$ is isogenous to a principally polarized abelian variety over $E$ (cf. [22], \S23, Corollary 1). Thus, based on the last two sentences, to prove Theorem 1.2 (and for the rest of the paper) we can assume that $A$ has a principal polarization $\lambda_A$.

\medskip\noindent
{\bf 2.3.1. Lemma.} {\it We assume that there exists a prime $p\in\dbN$ such that the group $H_{A\dbQ_p}^{\ad}$ is anisotropic i.e., its $\dbQ_p$-rank is $0$ (e.g., this holds if $H_A$ is a torus). Then there exists a finite field extension $E_1$ of $E$ such that $A_{E_1}$ extends to an abelian scheme over $O_{E_1}$.}

\medskip
\proof
The group $H^{\ad}_{A\dbQ_p}(\dbQ_p)$ is compact and thus it has no non-trivial unipotent element (see [27], Subsection 3.4). This implies that the group $H_A(\dbQ_p)$ has no non-trivial unipotent element. Thus the Lemma follows from [26], Theorem 1.6.1.\endproof

\medskip\noindent
{\bf 2.3.2. Remark.} From the classification of simple adjoint groups over $p$-adic fields (see [28], Table II, pp. 55--58) and the fact that each simple factor of $H_{A\dbC}^{\ad}$ is of classical Lie type (cf. [7], Table 2.3.8), one gets that the assumption that $H_{A\dbQ_p}^{\ad}$ has simple, anisotropic factors implies that $H_{A\dbC}^{\ad}$  has simple factors of $A_n$ Lie type for some $n\in\dbN$.
 
\bigskip\noindent
{\bf 2.4. Shimura pairs associated to $A$.} We recall that from now on $(A,\lambda_A)$ denotes a principally polarized abelian variety over a number field $E$. Let $L_A:=H_1(A(\dbC),\dbZ)$; it is a $\dbZ$-lattice of $W_A$. Let $\psi_A:L_A\otimes_{\dbZ} L_A\to\dbZ$ be the perfect, alternating form on $L_A$ induced by $\lambda_A$. If $W_A\otimes_{\dbQ} \dbC=F_A^{-1,0}\oplus F_A^{0,-1}$ is the Hodge decomposition defined by the homomorphism $h_A\colon\dbS\to \text{GL}_{W_A\otimes_{\dbQ} \dbR}$, let $\mu_A:\dbG_{m\dbC}\to \text{GL}_{W_A\otimes_{\dbQ} \dbC}$ be the Hodge cocharacter that fixes $F_A^{0,-1}$ and that acts via the identical character of $\dbG_{m\dbC}$ on $F_A^{-1,0}$.  We denote also by $h_A:\dbS\to H_{A\dbR}$ and $\mu_A:\dbG_{m\dbC}\to H_{A\dbC}$ the natural factorizations of $h_A$ and $\mu_A$ (respectively), cf. the very definition of $H_A$. Let $X_A$ be the $H_A(\dbR)$-conjugacy class of $h_A:\dbS\to H_{A\dbR}$. Let $S_A$ be the $\text{GSp}(W_A,\psi_A)(\dbR)$-conjugacy class of the homomorphism $\dbS\to\text{GSp}(W_A,\psi_A)_{\dbR}$ defined by $h_A$. It is well known that the pairs $(H_A,X_A)$ and $(\text{GSp}(W_A,\psi_A),S_A)$ are Shimura pairs and that we have an injective map $f_A:(H_A,X_A)\hookrightarrow (\text{GSp}(W_A,\psi_A),S_A)$ of Shimura pairs (see [6], [7], [16], and [18]). The Shimura variety $\Sh(\text{GSp}(W_A,\psi_A),S_A)$ is called a Siegel modular variety. A Shimura pair that admits an injective map into a Shimura pair that defines a Siegel modular variety, is called a Shimura pair of Hodge type; thus $(H_A,X_A)$ is a Shimura pair of Hodge type. Let $h_A^{\ad}:\dbS\to H_{A\dbR}^{\ad}$ be the composite of $h_A$ with the natural epimorphism $H_{A\dbR}\twoheadrightarrow H_{A\dbR}^{\ad}$. Let $X_A^{\ad}$ be the $H_A^{\ad}(\dbR)$-conjugacy class of $h_A^{\ad}$. The pair $(H_A^{\ad},X_A^{\ad})$ is called the adjoint Shimura pair of $(H_A,X_A)$. Similarly we define the adjoint Shimura pair $(G^{\ad},X^{\ad})$ of an arbitrary Shimura pair $(G,X)$. 

The $H_A(\dbC)$-conjugacy class $[\chi_{A\dbC}]$ of $\mu_A:\dbG_{m\dbC}\to H_{A\dbC}$ is defined over $\overline{\dbQ}$ and the Galois group $\Gal(\dbQ)$ acts on the corresponding $H_A(\overline{\dbQ})$-conjugacy class $[\chi_A]$ of cocharacters of $H_{A\overline{\dbQ}}$. The reflex field $E_A:=E(H_A,X_A)$ is the fixed field of the stabilizer subgroup of $[\chi_A]$ in $\Gal(\dbQ)$. Let $g:=\dim(A)$.

Let $N\Ge 3$ be an integer. Let $\psi_{A,N}:L_A/NL_A\otimes_{\dbZ/N\dbZ} L_A/NL_A\to \dbZ/N\dbZ$ be the reduction modulo $N$ of $\psi_A$. Let $(C,\lambda_C)$ be a principally polarized abelian scheme of relative dimension $g$ over a $\dbZ[{1\over N}]$-scheme $Y$. Let $\lambda_{C[N]}:C[N]\times_Y C[N]\to\mu_{NY}$ be the Weil pairing induced by $\lambda_C$. By a level-$N$ symplectic similitude structure of $(C,\lambda_C)$ we mean an isomorphism $\kappa:(L_A/NL_A)_Y\arrowsim C[N]$ of finite, \'etale group schemes over $Y$, such that there exists an element $\nu\in\mu_{NY}(Y)$ with the property that for all points $a$, $b\in (L_A/NL_A)_Y(Y)$ we have an identity $\nu^{\psi_{A,N}(a\otimes b)}=\lambda_{C[N]}(\kappa(a),\kappa(b))$ between elements of $\mu_{NY}(Y)$. 

Let $\scrA_{g,1,N}$ be the Mumford moduli scheme over $\dbZ[{1\over N}]$ that parameterizes principally polarized abelian schemes which are of relative dimension $g$ and which are equipped with a level-$N$ symplectic similitude structure, cf. [21], Theorems 7.9 and 7.10 naturally adapted to the case of level-$N$ symplectic similitude structures (instead of only level-$N$ structures). Let $(\scrA,\lambda_{\scrA})$ be the universal principally polarized abelian scheme over $\scrA_{g,1,N}$. 

Let $K(N):=\{h\in \text{GSp}(L_A,\psi_A)(\widehat\dbZ)|h\;\text{mod}\; N\;\text{is}\;\text{the}\;\text{identity}\}$. Let $K_A(N):=K(N)\cap H_A(\dbA_f)$. As $N\ge 3$, it is well-known that we can identify $\Sh(\text{GSp}(W_A,\psi_A),S_A)/K(N)={\scrA_{g,1,N}}_{\dbQ}$ (see [6], Proposition 4.17) and that the group $K(N)$ acts freely on $\Sh(\text{GSp}(W_A,\psi_A),S_A)$ (for instance, see [17], Subsections 2.10 to 2.14; this also follows from Serre Lemma of [22], Chapter IV, \S21, Theorem 5). To $f_A$ corresponds a finite morphism of $E_A$--schemes 
$$f_A(N):\Sh(H_A,X_A)/K_A(N)\to\Sh(\text{GSp}(W_A,\psi_A),S_A)_{E_A}/K(N)$$
which over $\dbC$ is obtained from the embedding $X_A\times H_A(\dbA_f)\hookrightarrow S_A\times \text{GSp}(W_A,\psi_A)(\dbA_f)$ between complex spaces via a natural passage to quotients, cf. [6], Corollary 5.4. As the group $K(N)$ acts freely on $\Sh(\text{GSp}(W_A,\psi_A),S_A)$, the group $K_A(N)$ also acts freely on $\Sh(H_A,X_A)$. This implies that the quotient epimorphism $\Sh(H_A,X_A)\twoheadrightarrow\Sh(H_A,X_A)/K_A(N)$ is a pro-\'etale cover and therefore $\Sh(H_A,X_A)/K_A(N)$ is a smooth $E_A$-scheme. 

Let $\scrN_N:=\scrN_{A,N}$ be the normalization of $(\scrA_{g,1,N})_{O_{E_A}[{1\over N}]}$ in the ring of fractions of $\Sh(H_A,X_A)/K_A(N)$; it is an $O_{E_A}[{1\over N}]$-scheme. Let $(\scrB,\lambda_{\scrB})$ be the pull back of $(\scrA,\lambda_{\scrA})$ to $\scrN_N$. 

\medskip\noindent
{\bf 2.4.1. On factors.} Let $\mu_A^{\ad}:\dbG_{m\dbC}\to H_{A\dbC}^{\ad}$ be the cocharacter naturally defined by $\mu_A$. If $H_t$ is a simple factor of $H_A^{\ad}$, let $\mu_t:\dbG_{m\overline{\dbQ}}\to H_{t\overline{\dbQ}}$ be a cocharacter whose extension to $\dbC$ is $H_t(\dbC)$-conjugate to the cocharacter of $H_{t\dbC}$ naturally defined by $\mu_A^{\ad}$. Let $h_{At}:\dbS\to H_{t\dbR}$ be the homomorphism naturally defined by $h_A$. Until Section 4, whenever the group $H_A^{\ad}$ is non-trivial we will denote by $H_0$ a fixed simple factor of $H_A^{\ad}$; therefore we will speak about the cocharacter $\mu_0:\dbG_{m\overline{\dbQ}}\to H_{0\overline{\dbQ}}$ and the homomorphism $h_{A0}:\dbS\to H_{0\dbR}$. 

\medskip\noindent
{\bf 2.4.2. On complex points.} We have $\Sh(H_A,X_A)(\dbC)=H_A(\dbQ)\backslash (X_A\times H_A(\dbA_f))$, cf. [7], Proposition 2.1.10 and Corollary 2.1.11. Let $u:=[x,h]\in \Sh(H_A,X_A)(\dbC)$, where $x\in X_A$ and $h\in H_A(\dbA_f)$. Let $W_A\otimes_{\dbQ} \dbC=F^{-1,0}_x\oplus F^{0,-1}_x$ be the Hodge decomposition defined by $x:\dbS\to H_{A\dbR}$ and let $L_h$ be the $\dbZ$-lattice of $W_A$ with the property that $h(L_A\otimes_{\dbZ} \widehat\dbZ)=L_h\otimes_{\dbZ} \widehat\dbZ$. We denote also by $u\in\Sh(H_A,X_A)/K_A(N)(\dbC)$ the image of $u$ through the epimorphism $\Sh(H_A,X_A)(\dbC)\twoheadrightarrow\Sh(H_A,X_A)/K_A(N)(\dbC)$ of sets. The complex torus associated to the abelian variety $B_u:=u^*(\scrB)$ is $F^{0,-1}_x\backslash (W_A\otimes_{\dbQ} \dbC)/L_h$, cf. Riemann's Theorem and the very construction of Siegel modular varieties (see [6], Theorem 4.7 and Example 4.16; in connection to $L_h$ see also [29], Subsection 4.1). The principal polarization $u^*(\lambda_{\scrB})$ of $B_u$ is uniquely determined by the property that it induces a perfect alternating form on $L_h$ which is a rational multiple of $\psi_A$ and which is a polarization of the Hodge $\dbQ$--structure on $W_A$ defined by $x$, cf. [29], Subsection 4.1. 

Let $C_A$ be the centralizer of $H_A$ in $\End(W_A)$. Due to Riemann's theorem we can naturally view $C_A$: (i) as a $\dbQ$--algebra of $\dbQ$--endomorphisms of any such pull pull back $B_u$ of $\scrB$, and (ii) as $\End(A_{\dbC})\otimes_{\dbZ} \dbQ$. We identify naturally $W_A=H_1(B_u(\dbC),\dbQ)$. Such an identification is unique up to isomorphisms $W_A\arrowsim W_A$ defined by elements of $H_A(\dbQ)$ and it is compatible with the natural actions of $C_A$.

\medskip\noindent
{\bf 2.4.3. Special points.} We now assume that $x\in X_A$ is a special point i.e., the homomorphism $x:\dbS\to H_{A\dbR}$ factors through the extension to $\dbR$ of a maximal torus $T_x$ of $H_A$. The Mumford--Tate group of $B_u$ is a reductive subgroup of $T_x$ and thus it is a torus. We have $u\in\im(\Sh(H_A,X_A)(\overline{\dbQ})\to \Sh(H_A,X_A)(\dbC))$, cf. [17], Theorem 1.8 applied to the special pair $(T_x,\{x\})$ of $(H_A,X_A)$. The set of all points $u=[x,h]\in\Sh(H_A,X_A)(\dbC)$ with $h\in H_A(\dbA_f)$, is Zariski dense in $\Sh(H_A,X_A)(\dbC)$ (cf. [6], Proposition 5.2).

\bigskip\noindent
{\bf 2.5. Toroidal compactifications.}  We consider a smooth, projective, toroidal compactification $\bar\scrA_{g,1,N}$ of $\scrA_{g,1,N}$ over $\dbZ[{1\over N}]$ such that the abelian scheme $\scrA$ over $\scrA_{g,1,N}$ extends to a semiabelian scheme $\bar\scrA$ over $\bar\scrA_{g,1,N}$ (cf. [10], Chapter IV, Theorem 6.7). We have:

\medskip
{\bf (i)} The fibres of $\bar\scrA$ over points of the complement of $\scrA_{g,1,N}$ in $\bar\scrA_{g,1,N}$, are semiabelian varieties that are not abelian varieties. 

\medskip
We consider the normalization $\bar\scrN_N:=\bar\scrN_{A,N}$ of $(\bar\scrA_{g,1,N})_{O_{E_A}[{1\over N}]}$ in the ring of fractions of $\Sh(H_A,X_A)/K_A(N)$; it is an $O_{E_A}[{1\over N}]$-scheme. As the morphism $f_A(N)$ is finite, as $\Sh(H_A,X_A)/K_A(N)$ is a normal (in fact even smooth) $E_A$-scheme, and as $\scrA_{g,1,N}$ is an open, Zariski dense subscheme of $\bar\scrA_{g,1,N}$, we have:

\medskip 
{\bf (ii)} The $O_{E_A}[{1\over N}]$-scheme $\scrN_N=\bar\scrN_N\times_{\bar\scrA_{g,1,N}} \scrA_{g,1,N}$ is an open, Zariski dense subscheme of $\bar\scrN_N$ and we have an identity $\scrN_{NE_A}=\Sh(H_A,X_A)/K_A(N)$.

\medskip
Let $\bar\scrB$ be the pull back of $\bar\scrA$ to $\bar\scrN_N$. Thus $\bar\scrB$ is a semiabelian scheme over $\bar\scrN_N$ whose restriction to $\scrN_N$ is the abelian scheme $\scrB$.

\bigskip\noindent
{\bf 2.6. Fact.} {\it The following two statements are equivalent:

\medskip
{\bf (a)} the $O_{E_A}[{1\over N}]$-scheme $\scrN_N$ is projective;

\smallskip
{\bf (b)} we have $\scrN_N=\bar\scrN_N$.

\medskip
Moreover, if these two statements hold, then there exists a finite field extension $E_1$ of $E$ such that $A_{E_1}$ extends to an abelian scheme over $O_{E_1}[{1\over N}]$.}

\medskip
\proof
As $O_{E_A}$ is an excellent ring (see [15], \S34), the scheme $\bar\scrN_N$ is a finite $(\bar\scrA_{g,1,N})_{O_{E_A}[{1\over N}]}$-scheme and therefore it is a projective $O_{E_A}[{1\over N}]$-scheme. But $\scrN_N$ is an open, Zariski dense subscheme of $\bar\scrN_N$ (cf. property 2.5 (ii)) and thus (a) is equivalent to (b).

To end the proof, it suffices to show that (a) implies the existence of a finite field extension $E_1$ of $E$ such that $A_{E_1}$ extends to an abelian scheme over $O_{E_1}[{1\over N}]$. Let $E_1$ be a number subfield of $\dbC$ that contains $i_E(E)$ and $E_A$ and such that there exists a level-$N$ symplectic similitude structure $\kappa$ of $(A,\lambda_A)_{E_1}$ whose pull back to $\dbC$ is defined by the canonical isomorphism $(L_A/NL_A)_{\dbC}\arrowsim A_{\dbC}[N]=({1\over N}L_A/L_A)_{\dbC}$. Let $v_A:\Spec(E_1)\to (\scrA_{g,1,N})_{O_{E_A}[{1\over N}]}$ be the morphism such that $(A,\lambda_A)_{E_1}=v_A^*((\scrA,\lambda_{\scrA})_{O_{E_A}[{1\over N}]})$ and the resulting level-$N$ symplectic similitude structure of $(A,\lambda_A)_{E_1}$ is $\kappa$. The composite of the morphism $\Spec(\dbC)\to\Spec(E_1)$ with $v_A$, is the complex point $[h_A,1_{W_A}]\in\text{Im}(f_A(N)(\dbC))$ (here $1_{W_A}$ is the identity element of $\text{GSp}(W_A,\psi_A)(\dbA_f)$). Thus, up to a replacement of $E_1$ by a finite field extension of it, $v_A$ factors through a morphism $u_A:\Spec(E_1)\to\scrN_N$. As (a) holds, from the valuative criterion of properness we get that $u_A$ extends to a morphism $u_{A,N}:\Spec(O_{E_1}[{1\over N}])\to\scrN_N$. The abelian scheme $u_{A,N}^*(\scrB)$  over $O_{E_1}[{1\over N}]$ extends $A_{E_1}$. \endproof

\bigskip\noindent
{\bf 2.7. Proposition.} {\it We assume that the principally polarized abelian scheme $(A,\lambda_A)$ over the number field $E$ is such that the $\dbQ$--rank of $H_A^{\ad}$ is $0$ (e.g., this holds if $A$ has compact factors). 

\medskip
{\bf (a)} Then the $E_A$-schemes $\scrN_{NE_A}$ and $\bar\scrN_{NE_A}$ coincide (i.e., we have $\scrN_{NE_A}=\bar\scrN_{NE_A}$).

\smallskip
{\bf (b)} We assume that $A$ has compact factors. We also assume that Theorem 1.2 holds for all abelian varieties over number fields which have compact factors. Then the $O_{E_A}[{1\over N}]$-scheme $\scrN_N$ is projective for every integer $N\Ge 3$.}

\medskip
\proof
As the $\dbQ$--rank of $H_A^{\ad}$ is $0$, the analytic space associated to $\Sh(H_A,X_A)_{\dbC}/K_A(N)$ is compact (see [2], Theorem 12.3 and Corollary 12.4). This implies that the analytic spaces associated to $\scrN_{N\dbC}$ and $\bar\scrN_{N\dbC}$ coincide. Thus $\scrN_{N\dbC}=\bar\scrN_{N\dbC}$. From this we get that (a) holds.

We prove (b). It suffices to show that the assumption that $\scrN_N\neq\bar\scrN_N$ leads to a contradiction, cf. Fact 2.6. As $\scrN_{NE_A}=\bar\scrN_{NE_A}$ and as $\scrN_N\neq\bar\scrN_N$, we get that $\bar\scrN_N$ has points with values in finite fields which do not belong to $\scrN_N$. As the morphism $\bar\scrN_N\to \Spec(O_{E_A}[{1\over N}])$ is flat, it has quasi-sections in the quasi-finite, flat topology of $\Spec(O_{E_A}[{1\over N}])$ whose images contain any a priori given point of $\bar\scrN_N$ with values in a finite field (cf. [12], Corollary (17.16.2)). From the last two sentences, we get that there exists a finite field extension $\tilde E$ of $E_A$ in $\dbC$ and a local ring $\tilde O$ of $O_{\tilde E}$ of mixed characteristic such that we have  a morphism $\tilde u:\Spec(\tilde O)\to\bar\scrN_N$ which does not factor through $\scrN_N$. Let $\tilde A$ be the generic fibre of $\tilde u^*(\bar\scrA)$; it is an abelian variety over $\tilde E$. To reach a contradiction, we can assume that $i_E(E)\subseteq \tilde E$. Let $H_{\tilde A}$ be the Mumford--Tate group of $\tilde A_{\dbC}$. 

To check that $\tilde A$ has compact factors, we can assume in this paragraph that the adjoint group $H_{\tilde A}^{\ad}$ is non-trivial. Let $H_{\tilde t}$ be an arbitrary simple factor of $H_{\tilde A}^{\ad}$. As the generic fibre of $\tilde u$ factors through $\Sh(H_A,X_A)/K_A(N)$, the group $H_{\tilde A}$ is the Mumford--Tate group defined by a homomorphism $\dbS\to\text{GL}_{W_A}$ which is an element of $X_A$. This implies that $\tilde H_A$ is naturally a subgroup of $H_A$. Thus we have natural inclusions $\Lie(H_{\tilde t})\subseteq\Lie(H_{\tilde A}^{\ad})\subseteq\Lie (H_A^{\ad})$. Let $H_t$ be a simple factor of $H_A^{\ad}$ with the property that the natural Lie homomorphism $\Lie(H_{\tilde t})\to\Lie(H_t)$ is a monomorphism of simple Lie algebras over $\dbQ$. As $A$ has compact factors, there exists a simple, compact factor $C_t$ of $H_{t\dbR}$. Let $C_{\tilde t}$ be a simple factor of $H_{\tilde t\dbR}$ such that the simple Lie algebra $\Lie(C_{\tilde t})$ over $\dbR$ is naturally a Lie subalgebra of $\Lie(C_t)$. The group $C_{\tilde t}$ is isogenous to a subgroup of the compact group $C_t$ and thus it is compact. This implies that the abelian variety $\tilde A$ has compact factors. 

Let $\tilde E_1$ be a finite field extension of $\tilde E$ such that $\tilde A_{\tilde E_1}$ extends to an abelian scheme over $O_{\tilde E_1}[{1\over N}]$, cf. our last hypothesis. Let $\tilde v_1$ be a prime of $\tilde E_1$ such that its local ring $\tilde O_1$ dominates $\tilde O$. As $\tilde A_{\tilde E_1}$ has good reduction with respect to $\tilde v_1$, the composite of the natural morphism $\Spec(\tilde O_1)\to\Spec(\tilde O)$ with $\tilde u$, factors through $\scrN_N$ (cf. property 2.5 (i)). Thus $\tilde u$ factors through $\scrN_N$. Contradiction. This proves (b).\endproof

\bigskip\noindent
{\bf 2.8. Lemma.} {\it Let $p\in\dbN$ be a prime that does not divide $N$. Let $k$ be an algebraic closure of the field $\dbF_p$ with $p$ elements. We assume that the $\dbQ$--rank of $H_A^{\ad}$ is $0$ (e.g., this holds if $A$ has compact factors). We also assume that there exists no morphism $q:\Spec(k[[x]])\to\bar\scrN_{N}$ with the property that it gives birth to morphisms $q_{\sp}:\Spec(k)\to\bar\scrN_{N}$ and $q_{\gen}:\Spec(k((x)))\to\bar\scrN_{N}$ that factor through $\bar\scrN_{N}\setminus\scrN_{N}$  and $\scrN_{N}$ (respectively). Then the complement $\bar\scrN_{N}\setminus\scrN_{N}$ has no points of characteristic $p$.}

\medskip
\proof
The only part of the proof of the Lemma which might be less well-known, is that $\bar\scrN_{N}\setminus\scrN_{N}$ does not contain the reduced scheme of any connected component of the special fibre in characteristic $p$ of $\scrN_{N}$. In the next two paragraphs we first check this property. 

We have $\scrN_{NE_A}=\bar\scrN_{NE_A}$, cf. Proposition 2.7 (a). The projective morphism $\bar\scrN_{N}\to\Spec(O_{E_{A}}[{1\over {N}}])$ is the composite of a projective morphism $n_{A,N}:\bar\scrN_{N}\to\Spec(O_{E_{A,N}}[{1\over {N}}])$ with connected fibres and of a finite morphism $\Spec(O_{E_{A,N}}[{1\over {N}}])\to\Spec(O_{E_{A}}[{1\over {N}}])$, cf. Stein's factorization theorem (see [13], Chapter III, Theorem 11.5); here $O_{E_{A,N}}$ is the ring of integers of a finite, \'etale $E_{A}$-algebra $E_{A,N}$. Let $\eta$ be an arbitrary point of $\Spec(O_{E_{A,N}}[{1\over {N}}])$ of characteristic $p$. Let $\scrF$ be the fibre of the morphism $n_{A,N}$ over $\eta$. 

Let $E_{\eta}$ be the field that is a direct factor of $E_{A,N}$ and such that $\eta$ is a point of the spectrum of the direct factor $O_{E_{\eta}}[{1\over {N}}]$ of $O_{E_{A,N}}[{1\over {N}}]$. Let $\scrE$ be the fibre of $n_{A,N}$ over $\Spec(E_{\eta})$; it is a connected component of $\scrN_{NE_{A}}=\bar\scrN_{NE_{A}}$. From the Zariski density part of Subsubsection 2.4.3, we get the existence of a finite field extension $E_{1\eta}$ of $E_{\eta}$ such that we have a morphism $u:\Spec(E_{1\eta})\to\scrN_{NE_{A}}=\bar\scrN_{NE_{A}}$ that factors through $\scrE$ and such that the Mumford--Tate group of a suitable (in fact of each) pull back of $u^*(\scrB_{E_{A}})$ to $\dbC$, is a torus. As $\bar\scrN_{N}$ is a projective $O_{E_{A}}[{1\over {N}}]$-scheme, the morphism $u$ extends to a morphism $\bar u:\Spec(O_{E_{1\eta}}[{1\over {N}}])\to\bar\scrN_{N}$. We can assume that the field $E_{1\eta}$ is such that the abelian variety $u^*(\scrB_{E_{A}})$ extends (cf. Lemma 2.3.1) to an abelian scheme over $O_{E_{1\eta}}[{1\over {N}}]$ which (cf. [10], Chapter I, Proposition 2.7) is the semiabelian scheme $\bar u^*(\bar\scrB)$. Thus $\bar u$ factors through $\scrN_{N}$, cf. property 2.5 (i). As $u$ factors through $\scrE$, $\text{Im}(\bar u)$ has a non-trivial intersection with $\scrF$. From the last two sentences, we get that the intersection $\scrF\cap\scrN_{N}$ is non-empty. 

But $\scrF\cap\scrN_{N}$ is an open subscheme of $\scrF$, cf. property 2.5 (ii). Our last hypothesis implies that the morphism $\scrF\cap \scrN_{N}\to\scrF$ is a closed embedding. As $\scrF$ is connected and has a non-empty intersection with $\scrN_{N}$, from the last two sentences we get that $\scrF=\scrF\cap \scrN_{N}$. Thus $\scrF$ is a closed subscheme of $\scrN_{N}$. Thus $\bar\scrN_{N}\setminus\scrN_{N}$ has no points of characteristic $p$.\endproof

\bigskip\smallskip
\noindent
{\boldsectionfont \S3. A construction}
\bigskip 

In this Section we assume that the abelian variety $A$ has compact factors (equivalently, that the Shimura pair $(H_A,X_A)$ has compact factors) and that the adjoint group $H_A^{\ad}$ is non-trivial. Let $(H_0,X_0)$ be a fixed simple factor of $(H_A^{\ad},X_A^{\ad})$; it has compact factors. The homomorphism $h_{A0}:\dbS\to H_{0\dbR}$ of Subsubsection 2.4.1 is an element of $X_0$ and in fact $X_0$ is the $H_0(\dbR)$-conjugacy class of $h_{A0}$. 

For another abelian variety $B$ over $E$, let $(H_B,X_B)$, $h_B:\dbS\to H_{B\dbR}$, $h_B^{\ad}:\dbS\to H^{\ad}_{B\dbR}$, and $(H_B^{\ad},X_B^{\ad})$ be the analogues of $(H_A,X_A)$, $h_A:\dbS\to H_{A\dbR}$, $h_A^{\ad}:\dbS\to H^{\ad}_{A\dbR}$, and $(H_A^{\ad},X_A^{\ad})$ (respectively) introduced in Subsection 2.4 but for $B$ instead of $A$. 

Let $p\in\dbN$ be a prime. In this Section we will prove the following result.

\bigskip\noindent
{\bf 3.1. Theorem.} {\it Up to a replacement of $E$ by a finite field extension of it, there exists a principally polarized abelian variety $(A_0,\lambda_{A_0})$ over $E$ such that the following two properties hold:

\medskip
{\bf (a)} we have an isomorphism $(H_{A_0}^{\ad},X_{A_0}^{\ad})\arrowsim (H_0,X_0)$ (to be viewed as an identity) with the property that the homomorphism $h_{A_0}^{\ad}:\dbS\to H_{A_0\dbR}^{\ad}=H_{0\dbR}$ is $h_{A0}$;

\smallskip
{\bf (b)} there exists an integer $N_0\Ge 3$ that is relatively prime to $p$ and such that we have $\scrN_{A_0,N_0}=\bar\scrN_{A_0,N_0}$ (i.e., and such that the statement 2.6 (b) holds for $(A_0,\lambda_{A_0},N_0)$).}

\medskip\noindent
{\bf 3.1.1. On the proof of Theorem 3.1.}
The proof of Theorem 3.1 is carried out in Subsections 3.2 to 3.4. The existence (up to a replacement of $E$ by a finite field extension of it) of a principally polarized abelian variety $(A_0,\lambda_{A_0})$ over $E$ such that the property 3.1 (a) holds, is an elementary consequence of [7], Proposition 2.3.10. The hard part is to show that we can choose $(A_0,\lambda_{A_0})$ so that the property 3.1 (b) holds as well. In order to achieve that the property 3.1 (b) holds, we will take $A_0$ so that the following two properties hold:
\medskip
{\bf (i)} the rank of the $\dbZ$-algebra $\End(A_{0\dbC})$ is sufficiently big;
\smallskip
{\bf (ii)} if $A_0$ extends to a semiabelian scheme over a local ring of $O_E$ of mixed characteristic $(0,p)$, then the natural action of $\End(A_{0\dbC})$ on the group of characters of the maximal torus of the special fibre of the semiabelian scheme extension, has some specific properties. 
\medskip
Due to properties (i) and (ii), the semiabelian scheme of the property (ii) will turn to be an abelian scheme. Condition 3.1 (b) will be implied by natural moduli analogues of properties (i) and (ii). In Subsection 3.2 we include notations that are essential for a review of the constructions of [7], Proposition 2.3.10 and for supplementing these constructions in order to be able to take $A_0$ so that the moduli analogues of properties (i) and (ii) hold. The mentioned review and supplementing process are the very essence of the construction of $A_0$ and are gathered in Lemma 3.3. In Subsection 3.4 we check that the moduli analogues of properties (i) and (ii) hold and we use this to end the proof of Theorem 3.1.
\bigskip\noindent
{\bf 3.2. Notations.} Let $F_0$ be a totally real number subfield of $\overline{\dbQ}\subseteq\dbC$ such that we have
$$H_0=\Res_{F_0/\dbQ} G_0,$$ 
with $G_0$ as an absolutely simple adjoint group over $F_0$ (cf. [7], Subsubsection 2.3.4); the field $F_0$ is unique up to $\Gal(\dbQ)$-conjugation. 
Let $i_{\text{nat}}:F_0\hookrightarrow \dbR$ be the embedding naturally defined by the inclusions $F_0\subseteq\overline{\dbQ}\subseteq\dbC$. We identify $\Hom(F_0,\dbR)=\Hom(F_0,\overline{\dbQ})$. 

Let $T_0$ be a maximal torus of $H_0$. Let $B_0$ be a Borel subgroup of $H_{0\overline{\dbQ}}$ that contains $T_{0\overline{\dbQ}}$. Let $\grD_0$ be the Dynkin diagram of $\Lie(H_{0\overline{\dbQ}})$ with respect to $T_{0\overline{\dbQ}}$ and $B_0$. We have $H_{0\overline{\dbQ}}=\prod_{i\in \Hom(F_0,\dbR)} G_0\times_{F_0} {}_i \overline{\dbQ}$. Thus $\grD_0$ is a disjoint union $\cup_{i\in \Hom(F_0,\dbR)} \grD_i$, where $\grD_i$ is the connected Dynkin diagram of $\Lie(G_0\times_{F_0} {}_i \overline{\dbQ})$ with respect to $(G_0\times_{F_0} {}_i \overline{\dbQ})\cap T_{0\overline{\dbQ}}$ and $(G_0\times_{F_0} {}_i \overline{\dbQ})\cap B_0$. Let $\grL_0$ be the Lie type of a (any) simple factor of $H_{0\dbC}$. As the group $H_{0\dbR}=\prod_{i\in \Hom(F_0,\dbR)} G_0\times_{F_0} {}_i \dbR$ has simple, compact factors and simple, non-compact factors (see Subsection 2.2), we have $[F_0:\dbQ]\Ge 2$. 

For a vertex $\alpha$ of $\grD_0$, let $\grg_{\alpha}$ be the 1 dimensional Lie subalgebra of $\Lie(B_0)$ that corresponds to $\alpha$. The Galois group $\Gal(\dbQ)$ acts on $\grD_0$ as follows. If $\gamma\in\Gal(\dbQ)$, then $\gamma(\alpha)$ is the vertex of $\grD_0$ defined by the identity $\grg_{\gamma(\alpha)}=i_{g_{\gamma}}(\gamma(\grg_{\alpha}))$, where $i_{g_{\gamma}}$ is the inner conjugation of $\Lie(H_{0\overline{\dbQ}})$ by an element $g_{\gamma}\in H_0(\overline{\dbQ})$ which normalizes $T_{0\overline{\dbQ}}$ and for which we have an identity $g_{\gamma}\gamma(B_0) g_{\gamma}^{-1}=B_0$. 

Let $\mu_0:\dbG_{m\overline{\dbQ}}\to H_{0\overline{\dbQ}}$ be as in Subsubsection 2.4.1. Let $\grV_0$ be the set of verticis of $\grD_0$ such that the unique cocharacter of $T_{0\overline{\dbQ}}$ that acts on $\grg_{\alpha}$ trivially if $\alpha\notin\grV_0$ and via the identical character of $\dbG_{m\overline{\dbQ}}$ if $\alpha\in\grV_0$, is $H_0(\overline{\dbQ})$-conjugate to $\mu_0$. Let $\grO_0$ be the set of verticis of $\grD_0$ formed by the orbit of $\grV_0$ under $\Gal(\dbQ)$. As $H_A$ is the smallest subgroup of $\text{GL}_{W_A}$ through which $h_A$ factors, the images of both $h_{A0}$ and $\mu_0$ are non-trivial. Thus the set $\grV_0$ is non-empty. The image of $h_{A0}$ in a simple factor $\scrF_0$ of $H_{0\dbR}$ is trivial if and only if the group $\scrF_0$ is compact (this is so as the centralizer of $\text{Im}(h_A^{\ad})$ in $H_{A\dbR}^{\ad}$ is a maximal compact subgroup of $H_{A\dbR}^{\ad}$, cf. [7], page 259). As $H_{0\dbR}$ has at least one simple, compact factor (cf. Definition 1.1), we get that: 

\medskip
-- there exists an element $i_0\in\Hom(F_0,\dbR)$ such that $\grV_0$ contains no vertex of $\grD_{i_0}$ (equivalently, such that the simple factor $G_0\times_{F_0} {}_{i_0} \dbR$ of $H_{0\dbR}$ is compact). 

\medskip
As the Hodge $\dbQ$--structure on $\Lie(H_0)$ defined by any element  $x_0\in X_0$ is of type $\{(-1,1),(0,0),(1,-1)\}$, for each $i\in\Hom(F_0,\dbR)$ the set $\grV_0$ contains at most one vertex of $\grD_i$. We know that $\grL_0$ is a classical Lie type, cf. [7], Table 2.3.8. Moreover, if $\grV_0$ contains a vertex of $\grD_i$, then with the standard notations of [4], Plates I to VI, this vertex is (cf. [7], Table 1.3.9): an arbitrary vertex if $\grL_0=A_n$, vertex $1$ if $\grL_0=B_n$, vertex $n$ if $\grL_0=C_n$, and an extremal vertex if $\grL_0=D_n$. The reflex field $E(H_0,X_0)$ of $(H_0,X_0)$ is the fixed field of the open subgroup of $\Gal(\dbQ)$ that stabilizes $\grV_0$, cf. [7], Proposition 2.3.6. 

If the Lie type $\grL_0$ is $A_n$, $B_n$, or $C_n$, then $(H_0,X_0)$ is said to be of $A_n$, $B_n$, or $C_n$ type. If $\grL_0=D_n$ with $n\Ge 5$, then $(H_0,X_0)$ is said to be:
\medskip
-- of $D_n^{\dbR}$ type, if for each embedding $i:F_0\hookrightarrow\dbR$, $\grO_0$ contains only the vertex $1$ of $\grD_i$;
\smallskip
-- of $D_n^{\dbH}$ type, if for each embedding $i:F_0\hookrightarrow\dbR$, $\grO_0$ contains the vertices $n-1$ or $n$ of $\grD_i$ but not the vertex $1$ of $\grD_i$. 
\medskip
If $\grL_0=D_4$, then $(H_0,X_0)$ is said to be of $D_4^{\dbR}$ (resp. of $D_4^{\dbH}$) if for each embedding $i:F_0\hookrightarrow\dbR$, $\grO_0$ contains only one (resp. exactly two) verticis of $\grD_i$; with the notations of [4], Plate IV, this vertex  (resp. these two verticis) will be chosen in what follows to be the vertex $1$ (resp. to be the last two verticis $3$ and $4$). 

The definition of the $A_n$, $B_n$, $C_n$, $D_n^{\dbH}$, and $D_n^{\dbR}$ types conforms with [7]. From [7], Table 2.3.8 we get that $(H_0,X_0)$ is of $A_n$, $B_n$, $C_n$, $D_n^{\dbH}$, or $D_n^{\dbR}$ type.

Let $\grS_0$ be the subset of verticis of $\grD_0$ defined as follows:

\medskip
-- if $(H_0,X_0)$ is of $A_n$ type, then $\grS_0$ is the set of all extremal verticis;

\smallskip
-- if $(H_0,X_0)$ is of $B_n$ (resp. $C_n$) type, then $\grS_0$ is the set of all verticis $n$ (resp. $1$);

\smallskip
-- if $(H_0,X_0)$ is of $D_n^{\dbH}$ (resp. $D_n^{\dbR}$) type, then $\grS_0$ is the set of all verticis $1$ (resp. $n-1$ and $n$).    

\medskip
The set $\grS_0$ is $\Gal(\dbQ)$-invariant (if $\grL_0=D_n$, this is implied by the very definitions of the $D_n^{\dbH}$ and $D_n^{\dbR}$ types). We identify the $\Gal(\dbQ)$-set $\grS_0$ with $\Hom(F_1,\dbC)$, where $F_1$ is an \'etale $F_0$-algebra of degree at most $2$. We have $[F_1:F_0]=2$ if and only if $(H_0,X_0)$ is either of $A_n$ type with $n\ge 2$ or of $D_n^{\dbR}$ type with $n\Ge 4$. The $F_0$-algebra $F_1$ is either a field of CM type (cf. [7], Subsubsection 2.3.4 (b) or the first paragraph of the proof of [7], Proposition 2.3.10) or a product of two fields isomorphic to $F_0$ and thus of CM type. If $[F_1:F_0]=2$ and $(H_0,X_0)$ is of $D_n^{\dbR}$ type with $n$ even, then $F_1$ is not necessarily a field. 

\medskip\noindent
{\bf 3.2.1. The field $K_0$.} Let $\overline{\dbQ_p}$ be an algebraic closure of $\dbQ_p$. We fix an identification between $\overline{\dbQ}$ and the algebraic closure of $\dbQ$ in $\overline{\dbQ_p}$ and we use it to identify naturally the set $\Hom(F_0,\dbR)=\Hom(F_0,\overline{\dbQ})$ with $\Hom(F_0,\overline{\dbQ_p})$. We write $F_0\otimes_{\dbQ} \dbQ_p=\prod_{j\in J} F_{0j}$ as a product of $p$-adic fields. Let $j_0\in J$ be the unique element such that to the embedding $i_0:F_0\hookrightarrow \dbR$ corresponds (under the identification $\Hom(F_0,\dbR)=\Hom(F_0,\overline{\dbQ_p})$) an embedding $F_0\hookrightarrow\overline{\dbQ_p}$ that factors through the composite embedding $F_0\hookrightarrow F_0\otimes_{\dbQ} \dbQ_p\twoheadrightarrow F_{0j_0}$. 

Let $v_{0j_0}$ be the prime of $F_0$ above $p$ such that the completion of $F_0$ with respect to $v_{0j_0}$ is the factor $F_{0j_0}$ of $F_0\otimes_{\dbQ} \dbQ_p$. Let $K_0$ be a totally imaginary quadratic extension of $F_0$ which is unramified above primes of $F_0$ that divide $p$ and which has only one prime $w_{0j_0}$ above $v_{0j_0}$. We have $[K_0:\dbQ]=2[F_0:\dbQ]\Ge 4$. Let $K_{0j_0}$ be the completion of $K_0$ with respect to $w_{0j_0}$; we have $[K_{0j_0}:F_{0j_0}]=2$. As $F_0$ is a totally real number field and as $K_0$ is a totally imaginary quadratic extension of $F_0$, the field $K_0$ is of CM type. 

\bigskip\noindent
{\bf 3.3. Lemma.} {\it We recall that $(H_0,X_0)$ is a simple factor of $(H_A^{\ad},X_A^{\ad})$. There exists a Shimura pair $(H_1,X_1)$ such that the following four properties hold:

\medskip
{\bf (i)} the adjoint Shimura pair of $(H_1,X_1)$ is $(H_0,X_0)$ and there exists $h_{A_0}\in X_1$ that maps naturally into the element $h_{A0}\in X_0$ introduced in Subsubsection 2.4.1;

\smallskip
{\bf (ii)} we have an injective map $f_1:(H_1,X_1)\hookrightarrow (\text{GSp}(W_1,\psi_1),S_1)$ into a Shimura pair that defines a Siegel modular variety;

\smallskip
{\bf (iii)} the torus $\scrT:=\Res_{K_0\otimes_{F_0} F_1/\dbQ} \dbG_{mK_0\otimes_{F_0} F_1}$ is naturally a subgroup of $\text{GL}_{W_1}$ that centralizes $H_1^{\der}$ and that makes $W_1$ to have a natural structure of a $K_0$-vector space;

\smallskip
{\bf (iv)} the torus $Z^0(H_1)$ is the torus of $\text{GL}_{W_1}$ generated by $Z(\text{GL}_{W_1})$ and by the maximal subtorus $\scrT_c$ of $\scrT$ which over $\dbR$ is compact.}

\medskip
\proof
The existence of the Shimura pair $(H_1,X_1)$ such that properties (i) and (ii) (resp. (iii) and (iv)) hold follows from the statement (resp. the proof) of [7], Proposition 2.3.10. We recall the details of the construction of $(H_1,X_1)$, in the form needed in what follows.
Let $Q_0$ be a maximal torus of $G_0$. Its rank is equal to the rank $n$ of $\grL_0$. To ease the notations, we can assume that we have an identity $T_0=\Res_{F_0/\dbQ} Q_0$ between tori of $H_0$. 

Let $E_0$ be the smallest subfield of $\dbC$ which contains $F_0$ and which has the property that the torus $Q_{0E_0}$ is split; it is a Galois extension of $F_0$ and the Galois group $\Gal(E_0/F_0)$ is a finite subgroup of $\text{GL}_{X^*(Q_{0E_0})}(\dbZ)$, where the group $X^*(Q_{0E_0})$ of characters of $Q_{0E_0}$ is viewed as a free $\dbZ$-module of rank $n$. 

We will consider a representation $\rho_0:G_{0E_0}^{\sc}\to \text{GL}_{V_0}$, where $V_0$ is an $E_0$-vector space of finite dimension. In what follows all the weights used are with respect to the maximal torus of $G_{0E_0}^{\sc}$ whose image in $G_{0E_0}$ is the maximal torus $Q_{0E_0}$ of $G_{0E_0}$. Depending on the type of $(H_0,X_0)$, we choose $\rho_0$ such that (cf. the definition of the subset $\grS_0$ of verticis of $\grD_0$) the following properties hold (see [4], Plates I to IV for the weights used):

\medskip
{\bf (v.a)} if $(H_0,X_0)$ is of $A_n$ type with $n\Ge 2$, then $\rho_0$ is the direct sum of the two faithful representations of $G_{0E_0}^{\sc}=\text{SL}_{n+1E_0}$ associated to the weights $\varpi_1$ and $\varpi_n$ (thus $\dim_{E_0}(V_0)=2n+2$); 

\smallskip
{\bf (v.b)} if $(H_0,X_0)$ is of $B_n$ type with $n\Ge 3$, then $\rho_0$ is the faithful spin representation of $G_{0E_0}^{\sc}=\text{Spin}_{2n+1E_0}$ associated to the weight $\varpi_n$ (thus $\dim_{E_0}(V_0)=2^n$);

\smallskip
{\bf (v.c)} if $(H_0,X_0)$ is of $C_n$ type with $n\Ge 1$, then $\rho_0$ is the faithful representation of $G_{0E_0}^{\sc}=\text{Sp}_{2nE_0}$ of dimension $2n$ associated to the weight $\varpi_1$ (thus $\dim_{E_0}(V_0)=2n$);

\smallskip
{\bf (v.d1)} if $(H_0,X_0)$ is of $D_n^{\dbH}$ type with $n\Ge 4$, then  $\rho_0$ is the representation of $G_{0E_0}^{\sc}$ of dimension $2n$ associated to the weight $\varpi_1$ (thus $\dim_{E_0}(V_0)=2n$);

\smallskip
{\bf (v.d2)} if $(H_0,X_0)$ is of $D_n^{\dbR}$ type with $n\Ge 4$, then  $\rho_0$ is the spin representation of $G_{0E_0}^{\sc}=\text{Spin}_{2nE_0}$ associated to the weights $\varpi_{n-1}$ and $\varpi_n$ (thus $\dim_{E_0}(V_0)=2^n$).

\medskip
Let $V_1$ be $V_0$ but viewed as a rational vector space; we keep in mind that $V_1$ has also a natural structure of an $E_0$-vector space and thus also of an $F_0$-vector space. As $H_0^{\sc}=\Res_{F_0/\dbQ} G^{\sc}_0$ is a subgroup of $\Res_{E_0/\dbQ} G^{\sc}_{0E_0}$, $V_1$ is naturally an $H_0^{\sc}$-module. Let $H_1^{\der}$ be the image of the natural representation $H_0^{\sc}\to \text{GL}_{V_1}$ over $\dbQ$; the adjoint group of $H_1^{\der}$ is $H_0$. The set of weights used in (v.a) to (v.d2) is stable under the natural action of $\Gal(F_0)$ on the abelian group of weights and as a $\Gal(F_0)$-set it can be identified with the $\Gal(F_0)$-set of verticis of $\grD_{i_{\text{nat}}}$ contained in $\grS_0$. This implies that the center of the double centralizer of $H_1^{\der}$ in $\text{GL}_{V_1}$ is the torus $\Res_{F_1/\dbQ} \dbG_{mF_1}$ (see Subsection 3.2 for $F_1$).

We take $W_1:=K_0\otimes_{F_0} V_1$ and we view it as a rational vector space. As $W_1$ has also a natural structure of a $K_0\otimes_{F_0} F_1$-module whose annihilator is trivial, the torus $\scrT:=\Res_{K_0\otimes_{F_0} F_1/\dbQ} \dbG_{mK_0\otimes_{F_0} F_1}$ is naturally a subgroup of $\text{GL}_{W_1}$. Moreover $W_1$ has a natural structure of a $K_0$-vector space. We will also identify $H_1^{\der}$ with a semisimple subgroup of $\text{GL}_{W_1}$ that commutes with $\scrT$. As $K_0$ and the simple factors of $F_1$ are fields of CM type (cf. Subsubsections 3.2 and 3.2.1), the maximal compact subtorus of $\scrT_{\dbR}$ is the extension to $\dbR$ of a subtorus $\scrT_c$ of $\scrT$. Let $Z^0(H_1)$ be the torus of $\text{GL}_{W_1}$ generated by $\scrT_c$ and $Z(\text{GL}_{W_1})$. The torus $Z^0(H_1)$ commutes with $H_1^{\der}$ and therefore there exists a unique reductive subgroup $H_1$ of $\text{GL}_{W_1}$ such that the notations match (i.e., the derived group of $H_1$ is $H_1^{\der}$ and the maximal torus of the center of $H_1$ is $Z^0(H_1)$). Thus the property (iv) holds. As $H_1$ commutes with $\scrT$, the property (iii) also holds. Let $H_2$ be the subgroup of $\text{GL}_{W_1}$ generated by $H_1^{\der}$ and $\scrT$; it contains $H_1$.

The existence of an injective map $f_1:(H_1,X_1)\hookrightarrow (\text{GSp}(W_1,\psi_1),S_1)$ such that the property (i) holds is part of the proof of [7], Proposition 2.3.10. We recall the part of loc. cit. that pertains to the existence of the element $h_{A_0}\in X_1$. We have $F_0\otimes_{\dbQ} \dbR=\prod_{i\in\Hom(F_0,\dbR)} \dbR$; for $i\in\Hom(F_0,\dbR)$ let $\pi_i$ be the idempotent of $F_0\otimes_{\dbQ} \dbR$ such that $\pi_i(F_0\otimes_{\dbQ} \dbR)$ is the factor $\dbR$ of $F_0\otimes_{\dbQ} \dbR$ that corresponds to $i$. Let $V(i):=\pi_iW_1\otimes_{\dbQ} \dbR$.  We have a direct sum decomposition $W_1\otimes_{\dbQ} \dbR=\oplus_{i\in\Hom(F_0,\dbR)} V(i)$ of $H_{1\dbR}^{\der}$-modules. We also have a direct sum decomposition 
$$W_1\otimes_{\dbQ} \dbC=\oplus_{i\in\Hom(F_0,\dbR)} W(i_1)\oplus W(i_2)\leqno (1)$$
of $H_{2\dbC}$-modules (and thus also of $H_{1\dbC}$-modules), where for each $i\in\Hom(F_0,\dbR)$ the elements $i_1$, $i_2\in\Hom(K_0,\dbC)$ extend $i\in\Hom(F_0,\dbC)$ and are listed in an a priori chosen order and where $V(i)\otimes_{\dbR} \dbC=W(i_1)\oplus W(i_2)$ is the natural decomposition into $K_0\otimes_{F_0} {}_i\dbR\otimes_{\dbR} \dbC$-modules. Each homomorphism $h_{A_0}:\dbS\to H_{2\dbR}$ normalizes $V(i)$ and thus gives birth to a homomorphism $h_{A_0,i}:\dbS\to \text{GL}_{V(i)}$. Moreover, each homomorphism $h_{A_0}:\dbS\to H_{2\dbR}$ that defines a Hodge $\dbQ$--structure on $W_1$ which has a (constant) weight, factors through $H_{1\dbR}$. 

We will choose a homomorphism $h_{A_0}:\dbS\to H_{2\dbR}$ such that the Hodge $\dbQ$--structure on $W_1$ is of type $\{(-1,0),(0,-1)\}$ (i.e., we have a natural Hodge decomposition $W_1\otimes_{\dbQ} \dbC=F_{A_0}^{-1,0}\oplus F_{A_0}^{0,-1}$ defined by $h_{A_0}$) and the following two additional properties hold:

\medskip
{\bf (vi.c)} if $i\in\Hom(F_0,\dbR)$ is such that $G_0\times_{F_0} {}_i \dbR$ is compact (for instance, if $i$ is $i_0$), then $h_{A_0,i}$ is fixed (i.e., centralized) by the image of $H_{2\dbR}$ in $\text{GL}_{V(i)}$, we have inclusions $W(i_1)\subseteq F_{A_0}^{-1,0}$ and $W(i_2)\subseteq F_{A_0}^{0,-1}$, and therefore $W(i_1)^*$ is included in the Hodge filtration $F^1_{A_0}(W_1^*\otimes_{\dbQ} \dbC)$ of $W_1^*\otimes_{\dbQ} \dbC$ defined by $h_{A_0}$; 

\smallskip
{\bf (vi.n)} if $i\in\Hom(F_0,\dbR)$ is such that $G_0\times_{F_0} {}_i \dbR$ is non-compact, then $h_{A_0,i}:\dbS\to \text{GL}_{V(i)}$ is the unique homomorphism such that the Hodge $\dbR$-structure on $V(i)$ is of type $\{(-1,0),(0,-1)\}$ and $h_{A_0,i}$ lifts the non-trivial homomorphism $\dbS\to G_0\times_{F_0} {}_i \dbR$ naturally defined by $h_{A0}$ (here $G_0\times_{F_0} {}_i \dbR$ is a simple factor of $H_{0\dbR}=\prod_{\tilde i\in\Hom(F_0,\dbR)} G_0\times_{F_0} {}_{\tilde i} \dbR$).

\medskip
See proof of [7], Proposition 2.3.10 for the explicit construction of $h_{A_0,i}$ of (vi.n); below we will only use (vi.c). We denote also by $h_{A_0}:\dbS\to H_{1\dbR}$ the factorization of $h_{A_0}$ through $H_{1\dbR}$ (the weight of the Hodge $\dbQ$--structure on $W_1$ defined by $h_{A_0}$ is $-1$). Let $X_1$ be the $H_1(\dbR)$-conjugacy class of $h_{A_0}:\dbS\to H_{1\dbR}$. From (vi.c) and (vi.n) we get that the property (i) holds. The existence of an injective map as in the property (ii) is a particular case of the argument for [7], Corollary 2.3.3. Thus the property (ii) also holds.\endproof 

\medskip\noindent
{\bf 3.3.1. Two extra properties.} In this Subsubsection we will use the notations of the proof of Lemma 3.3. From the property 3.3 (vi.c) we get that:

\medskip
{\bf (i)} for $i=i_0\in\Hom(F_0,\dbR)$ and each $x_1\in X_1$, $W(i_1)^*$ is included in the Hodge filtration $F^1_{x_1}(W_1^*\otimes_{\dbQ} \dbC)$ of $W_1^*\otimes_{\dbQ} \dbC$ defined by $x_1$.

\medskip
For each $i\in\Hom(F_0,\dbR)$, the real vector space $V(i)$ is the direct sum of all irreducible subrepresentations of the representation of $G_0^{\sc}\times_{F_0} {}_i \dbR$ on $W_1\otimes_{\dbQ} \dbR$. The group $\oplus_{\tilde i\in \Hom(F_0,\dbR)\setminus\{i\}} G_0^{\sc}\times_{F_0} {}_{\tilde i} \dbR$ fixes both $V(i)$ and $\psi_1$. Based on the last two sentences, the isomorphism $\delta_1:W_1\otimes_{\dbQ} \dbR\arrowsim W_1^*\otimes_{\dbQ} \dbR$ of $H_{1\dbR}^{\sc}$-modules (or $H_{1\dbR}^{\der}$-modules) naturally induced by $\psi_1$ has the property that for all $i\in\Hom(F_0,\dbR)$ it maps $V(i)$ onto the direct summand $V(i)^*$ of $W_1^*\otimes_{\dbQ} \dbR=(W_1\otimes_{\dbQ} \dbR)^*$ (i.e., we have $\delta_1(V(i))=V(i)^*$). Thus:

\medskip
{\bf (ii)} for all $i,\tilde i\in\Hom(F_0,\dbR)$ with $i\neq\tilde i$, the restriction of $\psi_1$ to $V(i)$ is non-degenerate and we have $\psi_1(V(i)\otimes_{\dbR} V(\tilde i))=0$.

\bigskip\noindent
{\bf 3.4. The proof of Theorem 3.1.} Let $H_{A_0}$ be the smallest subgroup of $H_1$ such that the homomorphism $h_{A_0}:\dbS\to H_{1\dbR}$ of the property 3.3 (i) factors through $H_{A_0\dbR}$. As either $2\pi i\psi_1$ or $-2\pi i\psi_1$ is a polarization of the Hodge $\dbQ$--structure on $W_1$ defined by $h_{A_0}$ (cf. property 3.3 (ii)), the group $H_{A_0}$ is reductive (cf. [8], Proposition 3.6). As $h_{A0}\in X_0$ is the image of $h_{A_0}\in X_1$, we have a natural identity $H_{A_0}^{\ad}=H_0$. Let $X_{A_0}$ be the $H_{A_0}(\dbR)$-conjugacy class of $h_{A_0}:\dbS\to H_{A_0\dbR}$. Let $L_{A_0}$ be a $\dbZ$-lattice of $W_1$ such that $\psi_1$ induces a perfect alternating form on $L_{A_0}$. 

We have an injective map $f_{A_0}:(H_{A_0},X_{A_0})\hookrightarrow (\text{GSp}(W_1,\psi_1),S_1)$ of Shimura pairs, cf. property 3.3 (ii). Let the 9-tuple $(E_{A_0},g_0,\scrA_{g_0,1,N_0},\scrN_{0N_0},\bar\scrN_{0N_0},\scrB_0,\lambda_{\scrB_0},\bar\scrB_0,K_{A_0}(N_0))$ 
be the analogue of the 9-tuple $(E_A,g,\scrA_{g,1,N},\scrN_N,\bar\scrN_N,\scrB,\lambda_{\scrB},\bar\scrB,K_A(N))$ formed by entries introduced in Subsections 2.4 and 2.5, but obtained in the context of the triple $(f_{A_0},L_{A_0},N_0)$ instead of the triple $(f_A,L_A,N)$; here the integer $N_0\Ge 3$ is relatively prime to $p$. Thus $E_{A_0}=E(H_{A_0},X_{A_0})$, $2g_0=\dim_{\dbQ}(W_1)$,  etc. From Subsection 2.4.2 applied in the context of the $7$-tuple $(L_{A_0},\scrA_{g_0,1,N_0},\scrN_{0N_0},\bar\scrN_{0N_0},\scrB_0,\lambda_{\scrB_0},K_{A_0}(N_0))$ instead of the $7$-tuple $(L_A,\scrA_{g,1,N},\scrN_{N},\bar\scrN_{N},\scrB,\lambda_{\scrB},K_{A}(N))$ and from the property 3.3 (iii), we get that we can naturally view $K_0$ as a $\dbQ$--algebra of $\dbQ$--endomorphisms of each pull pull back of the abelian scheme $\scrB_0$ via a $\dbC$-valued point of $\Sh(H_{A_0},X_{A_0})/K_{A_0}(N_0)$. This implies that, up to a replacement of $N_0$ by a positive integral power of it, we can view $K_0$ as a $\dbQ$--algebra of $\dbQ$--endomorphisms of the pull back of $\scrB_0$ to the spectrum of the ring of fractions of $\scrN_{0N_0}$ and thus also (cf. [10], Chapter I, Proposition 2.7) as a $\dbQ$--algebra of $\dbQ$--endomorphisms of either $\scrB_0$ or $\bar\scrB_0$. This represents the moduli analogue of the property 3.1.1 (i). 

The main point of the proof of Theorem 3.1 is to show that $\bar\scrN_{0N_0}\setminus\scrN_{0N_0}$ has no points of characteristic $p$ (the argument relies on Lemma 2.8 and it extends until Subsubsection 3.4.5). Let $k$ be an algebraic closure of the field $\dbF_p$.

\medskip\noindent
{\bf 3.4.1. An assumption.} We will show that the assumption that there exists a morphism $\break q:\Spec(k[[x]])\to\bar\scrN_{0N_0}$ with the property that it gives birth to morphisms $q_{\sp}:\Spec(k)\to\bar\scrN_{0N_0}$ and $q_{\gen}:\Spec(k((x)))\to\bar\scrN_{0N_0}$ that factor through $\bar\scrN_{0N_0}\setminus\scrN_{0N_0}$  and $\scrN_{0N_0}$ (respectively), leads to a contradiction (the argument will extend until Subsubsection 3.4.5). Let $C:=q^*(\bar\scrB_0)$; it is a semiabelian scheme over $k[[x]]$, whose generic fibre $C_{k((x))}$ is an abelian variety over $k((x))$ and whose special fibre $C_k$ is a semiabelian variety over $k$ that is not an abelian variety. Moreover, $\lambda_{C_{k((x))}}:=q_{\gen}^*(\lambda_{\scrB_0})$ is a principal polarization of $C_{k((x))}$. The next three Subsubsections represent the moduli analogue of the property 3.1.1 (ii). 

\medskip\noindent
{\bf 3.4.2. Notations.} Let $k_1:=\overline{k((x))}$. Let $B(k_1)$ be the field of fractions of the Witt ring $W(k_1)$ of $k_1$ and let $\sigma$ be its Frobenius automorphism. Let $(M,\phi,\psi_M)$ be the principally quasi-polarized $F$-isocrystal over $k_1$ of (the principally quasi-polarized $p$-divisible group of) $(C_{k((x))},\lambda_{C_{k((x))}})_{k_1}$. Thus $M$ is a $B(k_1)$-vector space of dimension $2g_0=\dim_{\dbQ}(W_1)$, $\phi:M\arrowsim M$ is a $\sigma$-linear automorphism, and $\psi_M:M\otimes_{B(k_1)} M\to B(k_1)$ is a non-degenerate alternating form which has the property that for all $a,b\in M$ we have an identity $\psi_M(\phi(a)\otimes\phi(b))=p\sigma(\psi_M(a\otimes b))$. Let $O$ be a finite discrete valuation ring extension of $W(k_1)$ such that we have a morphism $q_O:\Spec(O)\to\scrN_{0N_0}$ with the property that it gives birth to a morphism $q_{k_1}:\Spec(k_1)\to\scrN_{0N_0}$ which factors through $q_{\gen}$, cf. [12], Corollary (17.16.2) applied to the flat morphism $\scrN_{0N_0}\times_{\Spec(\dbZ)} \Spec(W(k_1))\to\Spec(W(k_1))$. By enlarging $O$, we can assume that the field of fractions $L$ of $O$ is naturally an algebra over the Galois extension of $\dbQ$ generated by $K_0$ and $E$; thus $K_0\otimes_{\dbQ} L\arrowsim L^{[K_0:\dbQ]}$ and therefore the set $\Hom(K_0,L)$ has $[K_0:\dbQ]$ elements. 

Let $(D,\lambda_D):=q_O^*(\scrB_0,\lambda_{\scrB_0})$; it is a principally polarized  abelian scheme over $O$. We fix an embedding $i_L:L\hookrightarrow\dbC$ that extends $i_E$; thus we can speak about $D_{\dbC}$. Viewing $K_0$ as a $\dbQ$--algebra of $\dbQ$--endomorphisms of the pull back $D$ of $\scrB_0$ (cf. Subsection 3.4), we get that we have a natural $\dbQ_p$-monomorphism $K_0\otimes_{\dbQ} \dbQ_p\hookrightarrow \End(M,\phi)$. Thus $M$ has a natural structure of a $K_0\otimes_{\dbQ} B(k_1)$-module and therefore also of an $F_0\otimes_{\dbQ} \dbQ_p$-module. As $F_0\otimes_{\dbQ} \dbQ_p=\prod_{j\in J} F_{0j}$, we have a unique decomposition of $F$-isocrystals over $k_1$
$$(M,\phi)=\oplus_{j\in J} (M_j,\phi)\leqno (2)$$
with the property that each $M_j$ is an $F_{0j}$-vector space.

We apply Subsection 2.4.2 in the context of $(L_{A_0},\scrA_{g_0,1,N_0},\scrN_{0N_0},\bar\scrN_{0N_0},\scrB_0,\lambda_{\scrB_0},K_{A_0}(N_0))$ instead of $(L_A,\scrA_{g,1,N},\scrN_{N},\bar\scrN_{N},\scrB,\lambda_{\scrB},K_{A}(N))$. Thus we have an identify $H_1(D({\dbC}),\dbQ)=W_1$ which is compatible with the natural $K_0$-actions. Moreover, the non-degenerate alternating form on $H_1(D({\dbC}),\dbQ)$ induced by $\lambda_D$ is a non-zero rational multiple of $\psi_1$.

\medskip\noindent
{\bf 3.4.3. Proposition.} {\it The $F$-isocrystal $(M_{j_0},\phi)$ has slopes $0$ and $1$ with multiplicity zero.} 

\medskip
\proof
Let $M^{00}_{j_0}$ be the $\dbQ_p$-vector subspace of $M_{j_0}$ formed by elements fixed by $\phi$. Thus $M^0_{j_0}:=M^{00}_{j_0}\otimes_{\dbQ_p} B(k_1)$ is the maximal $B(k_1)$-vector subspace of $M_{j_0}$ that is normalized by $\phi$ and such that all slopes of $(M^0_{j_0},\phi)$ are $0$. Obviously $M^{00}_{j_0}$ is a $K_0\otimes_{\dbQ} \dbQ_p$-module and thus $K_{0j_0}$ acts on $M^{00}_{j_0}$. As $M_{j_0}$ is an $F_{0j_0}$-vector space and as $K_{0j_0}$ is a field (see Subsubsection 3.2.1), $M^{00}_{j_0}$ is a $K_{0j_0}$-vector space.  

Let $F^1_L$ be the $L$-vector subspace of $M\otimes_{B(k_1)} L$ that defines the Hodge filtration of $M\otimes_{B(k_1)} L$ associated to the abelian variety $D_L$ via the functorial (in $D$) identification $M\otimes_{B(k_1)} L=H^1_{dR}(D_L/L)$ (see [5], Theorem 1.3). As $D$ is an abelian scheme over $O$, the triple $(M,\phi,F^1_L)$ is an admissible filtered module over $L$ in the sense of [11], Subsubsection 5.5.2 (cf. [11], Theorem of 6.1.4) and thus it is also a weakly-admissible filtered module over $L$ in the sense of [11], Definition 4.4.3 (cf. [11], Subsubsection 5.5.3). This implies that the Hodge polygon $\grp_{H}$ of $(M^0_{j_0},F^1_L\cap M^0_{j_0})$ is below the Newton polygon $\grp_N$ of $(M^0_{j_0},\phi)$, cf. [11], Proposition 4.4.2. As $\grp_N$ has all slopes $0$, we get that in fact $\grp_H=\grp_N$. Thus $(M^0_{j_0}\otimes_{B(k_1)} L)\cap F^1_L=0$. 

We fix an algebraic closure $\overline{L}$ of $L$ and we identify (to be compared with Subsubsection 3.2.1) $\overline{\dbQ}$ and $\overline{\dbQ_p}$ with their algebraic closures in $\overline{L}$. Thus we also identify the set $\Hom(F_0,\overline{\dbQ})=\Hom(F_0,\dbR)$ with the two sets $\Hom(F_0,\overline{\dbQ_p})$ and $\Hom(F_0,L)=\Hom(F_0,\overline{L})$ (resp. we identify the sets $\Hom(K_0,\overline{\dbQ})=\Hom(K_0,\dbC)$ and $\Hom(K_0,L)=\Hom(K_0,\overline{L})$).

Let $I_{0,p}^{(2)}$ (resp. $I_{0,p}$) be the subset of $\Hom(K_0,L)$ (resp. of $\Hom(F_0,L))$ formed by embeddings $K_0\hookrightarrow L$ (resp. $F_0\hookrightarrow L$) that have the property that (under them) the local ring $O$ of $L$ dominates the ring of integers of $F_{0j_0}$ in such a way that the resulting embedding $F_{0j_0}\hookrightarrow \overline{L}$ is defined by an embedding $F_{0j_0}\hookrightarrow\overline{\dbQ_p}$ which, up to $\Gal(\dbQ_p)$-conjugation, is (cf. Subsection 3.2) the element $i_0\in\Hom(F_0,\dbR)$. As $K_{0j_0}$ is a quadratic field extension of $F_{0j_0}$, the subset $I_{0,p}^{(2)}$ of $\Hom(K_0,L)$ has $[K_{0j_0}:\dbQ_p]$ elements and is $\Gal(K_0/F_0)$-invariant. Moreover the set $I_{0,p}$ has $[F_{0j_0}:\dbQ_p]$ elements and it is naturally identified with the quotient of $I_{0,p}^{(2)}$ under the action of $\Gal(K_0/F_0)$ on it. Let
$$M^{00}_{j_0}\otimes_{\dbQ_p} L=M^0_{j_0}\otimes_{B(k_1)} L=\oplus_{i_{0,L}\in \Hom(K_0,L)} M^{0i_{0,L}}_{j_0}$$ 
be the natural decomposition into $K_0\otimes_{\dbQ} L$-modules. As $M_{j_0}^{00}$ is a $K_{0j_0}$-vector space, each $M^{0i_{0,L}}_{j_0}$ with $i_{0,L}\in I_{0,p}^{(2)}$ is an $L$-vector space which is trivial if and only if $M^{00}_{j_0}=0$.

Formula (1) and the above two identifications $M\otimes_{B(k_1)} L=H^1_{dR}(D_L/L)$ and $H_1(D(\dbC),\dbQ)=W_1$ are functorial in $D$. We recall that $K_0$ is naturally a subfield of $\End(D)\otimes_{\dbZ} \dbQ$ (cf. Subsection 3.4) and that the principal polarization $\lambda_D$ of $D$ is defined by an isomorphism $D\arrowsim D^{\text{t}}$. From the last two sentences we get that we have natural identifications of $K_0$-vector spaces
$$M\otimes_{B(k_1)} \dbC=H^1_{dR}(D_{\dbC}/\dbC)=H^1(D({\dbC}),\dbQ)\otimes_{\dbQ} \dbC=W_1^*\otimes_{\dbQ} \dbC=\oplus_{i\in\Hom(F_0,\dbR)} W(i_1)^*\oplus W(i_2)^*,\leqno (3)$$ 
under which the following three properties hold (see [8], Section 1 for those properties that pertain to the relation between the de Rham and the Betti cohomologies of $D_{\dbC}$):

\medskip
{\bf (i)} $F^1_L\otimes_L \dbC$ gets identified with the Hodge filtration of $W_1^*\otimes_{\dbQ} \dbC$ defined by a point $x_1\in X_1$ that is naturally associated to $D_{\dbC}$;

\smallskip
{\bf (ii)} $M_{j_0}\otimes_{B(k_1)} \dbC$ gets identified with $(\oplus_{i\in I_{0,p}} V(i)^*)\otimes_{\dbR} \dbC$ (see proof of Lemma 3.3 for the $V(i)$'s);

\smallskip
{\bf (iii)} $\psi_M$ gets identified with a non-zero multiple of the non-degenerate alternating form on $W_1^*\otimes_{\dbQ} \dbC$ naturally induced by the non-degenerate alternating form $\psi_1$ on $W_1$.

\medskip
From properties (ii), (iii), and 3.3.1 (ii), we get that $\psi_M$ restricts to a non-degenerate alternating form on $M_{j_0}$ and thus it defines a principal quasi-polarization of $(M_{j_0},\phi)$. 
Therefore the $F$-isocrystal $(M_{j_0},\phi)$ has slopes $0$ and $1$ with equal multiplicities. Thus to end the proof of the Proposition, it suffices to show that the assumption that $(M_{j_0},\phi)$ has slope $0$ with positive multiplicity leads to a contradiction. We have $\dim_{B(k_1)}(M^0_{j_0})\Ge 1$; thus for $i_{0,L}\in I_{0,p}^{(2)}$ we have $M^{0i_{0,L}}_{j_0}\neq 0$. 

Under the identification $\Hom(F_0,\dbR)=\Hom(F_0,L)$, to the subset $I_{0,p}$ of $\Hom(F_0,L)$ corresponds a subset $I_0$ of $\Hom(F_0,\dbR)$ that contains the embedding $i_0:F_0\hookrightarrow\dbR$ of Subsection 3.2. Let $i_1,i_2:K_0\hookrightarrow\dbC$ be the two embeddings that extend $i:=i_0$ and that are listed in the same order as in the proof of Lemma 3.3. As the set $I_{0,p}^{(2)}$ is $\Gal(K_0/F_0)$-invariant, there exist elements $i_{1,p}$ and $i_{2,p}$ of $I_{0,p}^{(2)}$ that correspond to $i_1$ and $i_2$ (respectively) via the identification $\Hom(K_0,\dbC)=\Hom(K_0,L)$. From the property 3.3.1 (i) we get that the identifications of formula (3) give birth to inclusions 
$$M^{0i_{1,p}}_{j_0}\otimes_{B(k_1)} \dbC\subseteq W(i_1)^*\subseteq F^1_L\otimes_L \dbC.$$ 
\indent
Thus $M^{0i_{1,p}}_{j_0}\otimes_{B(k_1)} L\subseteq F^1_L$. Therefore $(M^0_{j_0}\otimes_{B(k_1)} L)\cap F^1_L\supseteq M^{0i_{1,p}}_{j_0}\otimes_{B(k_1)} L\supsetneqq 0$. This contradicts the identity $(M^0_{j_0}\otimes_{B(k_1)} L)\cap F^1_L=0$.\endproof 

\medskip\noindent
{\bf 3.4.4. The study of $C$.} We now use Proposition 3.4.3 to reach the contradiction promised in Subsubsection 3.4.1. Let $T_k$ be the maximal torus of the semiabelian variety $C_k$. As $C_k$ is not an abelian scheme, we have $1\Le\dim(T_k)$. Let $K_{0\dbZ}:=K_0\cap\End(C)$ (the intersection being taken inside $\End(C)\otimes_{\dbZ} \dbQ$); it is a $\dbZ$-order of $K_0$. As $K_{0\dbZ}$ acts on $C$, it also acts on $T_k$ and thus also on the free $\dbZ$-module $X^*(T_k)$ of characters of $T_k$. Let $m,l\in\dbN$. 

There exists a unique torus $T_{k,l}$ of $C_{k[[x]]/(x^l)}$ which lifts $T_k$, cf. [9], Exp. IX, Theorem 3.6 bis. Loc. cit. implies that we have a canonical identification $T_{k,l}=T_k\times_k {k[[x]]/(x^l)}$ that lifts the identity automorphism of $T_k$. Thus $T_{k,l}[p^m]=T_k[p^m]_{k[[x]]/(x^l)}$ is naturally a closed subgroup scheme of $C_{k[[x]]/(x^l)}[p^m]$. Due to the uniqueness property of $T_{k,l}$, the torus $T_{k,l+1}$ lifts $T_{k,l}$. Thus by passing to the limit $l\to\infty$, we get that $T_k[p^m]_{k[[x]]}$ is naturally identified with a closed subgroup scheme of $C_{k[[x]]}[p^m]$ and thus that $T_k[p^m]_{k((x))}$ is naturally identified with a closed subgroup scheme of $C_{k((x))}[p^m]$. To check that these last identifications are functorial, it suffices to show that for each closed, semiabelian subscheme $C^\prime$ of $C^2$, the unique subtorus $T^\prime_{k,l}$ of $C^\prime_{k[[x]]/(x^l)}$ that lifts the maximal torus $T^\prime_k$ of $C^\prime_k$, is a subtorus of $T^2_{k,l}$. As $T^\prime_k$ is a subtorus of $T^2_k$, from the uniqueness part of loc. cit. we get that: (i) there exists a unique subtorus $T^{\prime\prime}_{k,l}$ of $T^2_{k,l}$ that lifts $T^\prime_k$, and (ii) we have an identity $T^\prime_{k,l}=T^{\prime\prime}_{k,l}$ of subtori of $C^2_{k[[x]]/(x^l)}$. Thus $T^\prime_{k,l}$ is a subtorus of $T^2_{k,l}$.

The closed embedding homomorphism $\Theta_m:T_k[p^m]_{k((x))}\hookrightarrow C_{k((x))}[p^m]$ is compatible with the natural $K_{0\dbZ}$-actions, cf. the functorial part of the previous paragraph. By taking $m\to\infty$ we get that we have a monomorphism $\Theta_{\infty}:T_k[p^{\infty}]_{k((x))}\hookrightarrow C_{k((x))}[p^{\infty}]$ of $p$-divisible groups over $k((x))$ which is compatible with the $K_{0\dbZ}$-actions. The $F$-isocrystal of the $p$-divisible group $T_k[p^{\infty}]_{k_1}$ is the pair $(X^*(T_k)\otimes_{\dbZ} B(k_1),1_{X^*(T_k)}\otimes p\sigma)$ (such an identification is unique up to a scalar multiplication by a unit of $\dbZ_p$). To $\Theta_{\infty}$ corresponds an epimorphism of $F$-isocrystals over $k_1$
$$\theta_{\infty}:M\twoheadrightarrow X^*(T_k)\otimes_{\dbZ} B(k_1)\leqno (4)$$  
which is compatible with the $K_{0\dbZ}$-actions. As $M$ is a $K_0$-vector space (cf. formula (3)), $K_{0\dbZ}$ can not act trivially on a quotient of $M$ of positive dimension. Due to this and the existence of the epimorphism $\theta_{\infty}$ (see (4)), the action of $K_{0\dbZ}$ on $X^*(T_k)$ is non-trivial i.e., it is defined by a $\dbZ$-monomorphism $K_{0\dbZ}\hookrightarrow \End(X^*(T_k))$. Due to this property, the unique decomposition $X^*(T_k)\otimes_{\dbZ} \dbQ_p=\oplus_{j\in J} X^*(T_k)_j$ with the property that each $X^*(T_k)_j$ is an $F_{0j}$-vector space (to be compared with formula (2)), is such that every $F_{0j}$-vector space $X^*(T_k)_j$ is non-zero. From this and the existence of the epimorphism $\theta_{\infty}$ (see (4)), we get that for the element $j_0\in J$ we have an epimorphism
$$\theta_{\infty,j_0}:(M_{j_0},\phi)\twoheadrightarrow (X^*(T_k)_{j_0}\otimes_{\dbQ_p} B(k_1),1_{X^*(T_k)_{j_0}}\otimes p\sigma)$$
of $F$-isocrystals over $k_1$. Thus $(M_{j_0},\phi)$ has slope $1$ with positive multiplicity. This contradicts Proposition 3.4.3 i.e., the assumption of Subsubsection 3.4.1 leads to a contradiction. In other words, a morphism $q:\Spec(k[[x]])\to\bar\scrN_{0N_0}$ as in Subsubsection 3.4.1 does not exist. Thus the complement $\bar\scrN_{0N_0}\setminus\scrN_{0N_0}$ has no points of characteristic $p$, cf. Lemma 2.8 (applied to $(A_0,\lambda_{A_0})$) and the fact that $A_0$ (equivalently, $(H_0,X_0)$) has compact factors. 

\medskip\noindent
{\bf 3.4.5. End of the proof of Theorem 3.1.} We recall that $\bar\scrN_{0N_0}$ is a projective $O_{E_{A_0}}[{1\over N_0}]$-scheme and that (cf. Proposition 2.7 (a) applied in the context of $f_{A_0}$) we have $\bar\scrN_{0N_0E_{A_0}}=\scrN_{0N_0E_{A_0}}$. From this and the identity $(\bar\scrN_{0N_0}\setminus\scrN_{0N_0})_{\text{red}\dbF_p}=\emptyset$ (cf. end of Subsection 3.4.4), we get that, by replacing $N_0$ with $N_0c_0$ for some number $c_0\in\dbN$ relatively prime to $p$, we can assume that in fact we have $\scrN_{0N_0}=\bar\scrN_{0N_0}$. 

By replacing $E$ with a finite field extension of it, we can assume (see proof of Fact 2.6) that there exists a morphism $u_A:\Spec(E)\to\scrN_N$ such that $(A,\lambda_A)=u_A^*(\scrB,\lambda_{\scrB})$ and:

\medskip
{\bf (i)} the composite of the morphism $\Spec(\dbC)\to\Spec(E)$ defined by $i_E$ with $u_A$ is the point of $\scrN_N(\dbC)=\Sh(H_A,X_A)/K_A(N)(\dbC)=H_A(\dbQ)\backslash (X_A\times H_A(\dbA_f)/K_A(N))$ defined by the equivalence class $[h_A,1_{W_A}]$ (here $1_{W_A}$ is the identity element of $H_{A}(\dbA_f)$). 

\medskip
By replacing $N$ and $N_0$ with $Nc$ and $N_0c_0$, where $c$ and $c_0$ are natural numbers prime to $p$, we can assume that there exists a compact open subgroup $K_0$ of $H_0(\dbA_f)$ such that the images of both $K_A(N)$ and $K_{A_0}(N_0)$ in $H_0(\dbA_f)$, are contained in $K_0$. We have functorial morphisms $\Sh(H_A,X_A)/K_A(N)\to\Sh(H_0,X_0)/K_0$ and $\Sh(H_{A_0},X_{A_0})/K_{A_0}(N_0)\to\Sh(H_0,X_0)/K_0$, the last one being finite. Based on this and the property (i), by replacing $E$ with a finite field extension of it, we can assume that there exists a morphism $u_{A_0}:\Spec(E)\to\scrN_{0N_0}$ such that the $E$-valued points of $\Sh(H_0,X_0)/K_0$ naturally defined by $u_A$ and $u_{A_0}$ coincide and moreover:

\medskip
{\bf (ii)} the composite of the morphism  $\Spec(\dbC)\to\Spec(E)$ defined by $i_E$ with $u_{A_0}$ is the point of $\scrN_{0N_0}(\dbC)=\Sh(H_{A_0},X_{A_0})/K_{A_0}(N_0)(\dbC)=H_{A_0}(\dbQ)\backslash (X_{A_0}\times H_{A_0}(\dbA_f)/K_{A_0}(N_0))$ defined by the equivalence class $[h_{A_0},1_{W_1}]$ (here $1_{W_1}$ is the identity element of $H_{A_0}(\dbA_f)$). 

\medskip
Let $(A_0,\lambda_{A_0}):=u_{A_0}^*(\scrB_0,\lambda_{\scrB_0})$. We can naturally identify the triple $(W_1,\psi_1,X_1)$ with $(W_{A_0},\psi_{A_0},X_{A_0})$ and thus the notations for $(A_0,\lambda_{A_0})$ and for the following $9$-tuple $(E_{A_0},g_0,\scrA_{g_0,1,N_0},\scrN_{0N_0},\bar\scrN_{0N_0},\scrB_0,\lambda_{\scrB_0},\bar\scrB_0,K_{A_0}(N_0))$ match. Based on the property (ii) and the definition of $H_{A_0}$ in the beginning of Subsection 3.4, we get that the Mumford--Tate group of $A_{0\dbC}$ is $H_{A_0}$. Thus, as $h_{A_0}$ lifts $h_{A0}$, property 3.1 (a) holds. As $\scrN_{0N_0}=\bar\scrN_{0N_0}$, property 3.1 (b) also holds. This ends the proof of Theorem 3.1. As $\scrN_{0N_0}=\bar\scrN_{0N_0}$, property 3.1.1 (ii) holds trivially. As we have a natural monomorphism $K_0\hookrightarrow  \End(A_{0\dbC})$ (cf. Subsection 3.4) and as $[K_0:\dbQ]=2[F_0:\dbQ]\Ge 4$, the property 3.1.1 (i) also holds.\endproof  

\bigskip\smallskip
\noindent
{\boldsectionfont \S4. Proof of the Theorem 1.2, examples, and applications}
\bigskip

In Subsection 4.1 we prove Theorem 1.2. Example 4.2 is completely new. Corollary 4.3 is an equivalent form of Theorem 1.2. In Subsections 4.4 to 4.6 we apply Corollary 4.3 to N\'eron models and to integral models of Shimura varieties of preabelian type. We use the notations of the first paragraph of Section 1 and of Subsection 2.4.

\bigskip\noindent
{\bf 4.1. Proof of Theorem 1.2.} In this Subsection we use the embedding $i_E:E\hookrightarrow \dbC$ to view $E$ as a subfield of $\dbC$. Accordingly, all finite field extensions of $E$ will be viewed as subfields of $\dbC$ that contain $E$ and their composites will be taken inside $\dbC$. To prove Theorem 1.2, we can assume that the abelian variety $A$ has a principal polarization $\lambda_A$ (cf. Subsection 2.3) and that the group $H_A^{\ad}$ is non-trivial (cf. Lemma 2.3.1). Let $N_A\in\dbN$ be such that $A$ extends to an abelian scheme over $O_E[{1\over {N_A}}]$. Suppose that for each prime  divisor $p$ of $N_A$, there exists a finite field extension $E_p$ of $E$ such that $A_{E_p}$ has good reduction with respect to all primes of $E_p$ that divide $p$. If $E_1$ is the composite field of $E$ and of all the fields $E_p$'s with $p$ a prime divisor of $N_A$, then $A_{E_1}$ extends to an abelian scheme over $O_{E_1}$. Thus to end the proof of Theorem 1.2, we only need to show that the finite field extension $E_p$ of $E$ exists for all prime divisors $p$ of $N_A$. 

To check this, we can replace $E$ by any finite field extension of it. Let 
$$(H_A^{\ad},X_A^{\ad})=\prod_{t\in\grT} (H_t,X_t)$$ 
be the product decomposition into simple, adjoint Shimura pairs. By replacing $E$ with a finite field extension of it, based on Theorem 3.1 we can assume that for each $t\in\grT$ there exists a principally polarized abelian variety $(A_t,\lambda_{A_t})$ over $E$ such that:

\medskip
{\bf (i)} we have an identity $(H_{A_t}^{\ad},X_{A_t}^{\ad})=(H_t,X_t)$ with the property that the homomorphism $h_{A_t}^{\ad}:\dbS\to H_{A_t\dbR}^{\ad}=H_{t\dbR}$ is the homomorphism $h_{At}$ of Subsubsection 2.4.1, and 

\smallskip
{\bf (ii)} there exists an integer $N_t\Ge 3$ that is relatively prime to $p$ and such that we have $\scrN_{A_t,N_t}=\bar\scrN_{A_t,N_t}$ (i.e., and such that the statement 2.6 (b) holds for $(A_t,\lambda_{A_t},N_t)$).

\medskip
Let $E_{p,t}$ be a finite field extension of $E$ such that $A_{tE_{p,t}}$ has good reduction with respect to all primes of $E_{p,t}$ that divide $p$, cf. property (ii) and the last part of Fact 2.6 applied to $A_t$. Let $\tilde E_p$ be the composite field of $E_{p,t}$'s, with $t\in\grT$. Let $B:=\prod_{t\in\grT} A_t$; it is an abelian variety over $E$ with the property that $B_{\tilde E_p}$ has good reduction with respect to all primes of $\tilde E_p$ that divide $p$. The group $H_A^{\ad}$ is the smallest subgroup of $H_A^{\ad}=\prod_{t\in\grT} H_t$ with the property that $h_A^{\ad}=\prod_{t\in\grT} h_{At}$ factors through $H^{\ad}_{A\dbR}$, cf. the very definition of $H_A$. The Mumford--Tate group $H_B$ is a subgroup of $\prod_{t\in\grT} H_{A_t}$ that surjects onto all groups $H_{A_t}$. This implies that $H_B^{\ad}$ is the smallest subgroup of $\prod_{t\in\grT} H_{A_t}^{\ad}$ with the property that $h_B^{\ad}=\prod_{t\in\grT} h_{A_t}^{\ad}$ factors through $H^{\ad}_{B\dbR}$. From the last two sentences and the property (i) we get that:

\medskip
{\bf (iii)} we have identifications $(H_{B}^{\ad},X_{B}^{\ad})=\prod_{t\in\grT} (H_t,X_t)=(H_A^{\ad},X_A^{\ad})$ with the property that the homomorphism $h_{B}^{\ad}:\dbS\to H_{B\dbR}^{\ad}=H_{A\dbR}^{\ad}$ is the homomorphism $h_{A}^{\ad}$ of Subsection 2.4.

\medskip
The reductive group $H_{A\times_E B}$ is a subgroup of $H_A\times_{\dbQ} H_B$ whose adjoint is (cf. property (iii)) isomorphic to $H_A^{\ad}=H_B^{\ad}$. As $B_{\tilde E_p}$ has good reduction with respect to all primes of $\tilde E_p$ that divide $p$, from the property (iii) and [26], Proposition 4.1.2 we get that there exists a finite field extension $E_p$ of $\tilde E_p$ such that $A_{E_p}$ has good reduction with respect to all primes of $E_p$ that divide $p$. Thus the finite field extension $E_p$ of $E$ exists for each prime divisor $p$ of $N_A$. This ends the proof of Theorem 1.2.\endproof

\bigskip\noindent
{\bf 4.2. Example.} Let $F$ be a totally real, cubic, Galois extension of $\dbQ$; for instance, we can take $F$ to be $\dbQ(\zeta_7+\zeta_7^{-1})$, where $\zeta_7$ is a primitive root of $1$ of order $7$. We assume that $H_A^{\ad}$ is a simple group of the form $\Res_{F/\dbQ} G$ for some absolutely simple, adjoint group $G$ over $F$ of $B_n$ (resp. of $D_n$) Dynkin type with $n\Ge 2$ (resp. with $n\Ge 4$). We also assume that the product decomposition $H_{A\dbR}^{\ad}=\scrF_1\times_{\dbR} \scrF_2\times_{\dbR} \scrF_3$ into simple factors is such that $\scrF_1$ and $\scrF_2$ are non-compact and $\scrF_3$ is compact. Thus $A$ has compact factors and therefore the Morita conjecture holds for $A$, cf. Theorem 1.2.

We identify $\Lie(H_{A\dbC}^{\der})$ with $\Lie(H_{A\dbC}^{\ad})=\Lie(\scrF_{1\dbC})\oplus \Lie(\scrF_{2\dbC})\oplus\Lie(\scrF_{3\dbC})$. We check that the representation of $\Lie(H_{A\dbC}^{\ad})$ on $W_A\otimes_{\dbQ} \dbC$ is free of tensor products i.e., it is a direct sum of irreducible representations of either $\Lie(\scrF_{1\dbC})$ or $\Lie(\scrF_{2\dbC})$ or $\Lie(\scrF_{3\dbC})$. We show that the assumption that this is not true, leads to a contradiction. This assumption implies that there exists $s\in\{2,3\}$ and $\scrL_1,\ldots,\scrL_s\in\{\Lie(\scrF_{1\dbC}),\Lie(\scrF_{2\dbC}),\Lie(\scrF_{3\dbC})\}$ such that a suitable simple $\Lie(H_{A\dbC}^{\ad})$-submodule $\scrW_0$ of $W_A\otimes_{\dbQ} \dbC$ is a tensor product $\scrW_1\otimes_{\dbC}\cdots\otimes_{\dbC} \scrW_s$, where $\scrW_r$ is a simple $\scrL_r$-module for all $r\in\{1,\ldots,s\}$. As the representation of $\Lie(H_{A\dbC}^{\ad})$ on $W_A\otimes_{\dbQ} \dbC$ is defined over $\dbQ$ and as $F$ is a cubic Galois extension of $\dbQ$, we can choose $\scrW_0$ such that we have $\scrL_1=\Lie(\scrF_{1\dbC})$ and $\scrL_2=\Lie(\scrF_{2\dbC})$. As the images of $\mu_A$ in $\scrF_{1\dbC}$ and $\scrF_{2\dbC}$ are non-trivial,  for $s=2$ (resp. for $s=3$) the Hodge filtration $(F^a(\scrW_0))_{a\in\dbZ}$ of $\scrW_0$ defined by $\mu_A$ is the tensor product of non-trivial Hodge filtrations of $\scrW_1$ and $\scrW_2$ (resp. is the tensor product of non-trivial Hodge filtrations of $\scrW_1$ and $\scrW_2$ and of a trivial Hodge filtration of $\scrW_3$); here by a trivial filtration of $W_r$ we mean a filtration of $W_r$ that does not contain any proper subspace of $W_r$. We easily get that there exists $a\in\dbZ\setminus\{-1,0\}$ such that $F^a(\scrW_0)/F^{a+1}(\scrW_0)\neq 0$. Thus $\mu_A$ does not act on $F^a(\scrW_0)/F^{a+1}(\scrW_0)$ via either the trivial or the identical character of $\dbG_{m\dbC}$. This contradicts the very definition of $\mu_A$.  

As the representation of $\Lie(H_{A\dbC}^{\ad})$ on $W_A\otimes_{\dbQ} \dbC$ is free of tensor products, from Remark 2.3.2 and [26], Proposition 2.2.3 we get that the results of [26] (which pertain to perfectly tens-twisted representations defined in [26], Definition 2.2.2) do not imply that the Morita conjecture holds for $A$. Thus our example is completely new. 

\medskip
Based on Proposition 2.7 (b), we have the following equivalent form of Theorem 1.2.

\bigskip\noindent
{\bf 4.3. Basic Corollary.} {\it We assume that the principally polarized abelian scheme $(A,\lambda_A)$ over the number field $E$ is such that $A$ has compact factors. Let $E_A$, $g$, $\scrA_{g,1,N}$, and $\scrN_N$ be as in Subsection 2.4. Then the normal $O_{E_A}[{1\over N}]$-scheme $\scrN_N$ is projective and it is a finite scheme over $(\scrA_{g,1,N})_{O_{E_A}[{1\over N}]}$.}

\bigskip\noindent
{\bf 4.4. N\'eron models.} Let $\scrK$ be the field of fractions of an integral Dedekind ring $\scrD$. Let $Z_{\scrK}$ be a smooth, separated $\scrK$-scheme of finite type. We recall (cf. [3], page 12) that a N\'eron model of $Z_{\scrK}$ over $\scrD$ is a smooth, separated $\scrD$-scheme $Z$ of finite type that has $Z_{\scrK}$ as its generic fibre and that satisfies the following universal (N\'eron mapping) property: 

\medskip
{\it for each smooth $\scrD$-scheme $Y$ and each $\scrK$-morphism $y_{\scrK}:Y_{\scrK}\to Z_{\scrK}$, there exists a unique morphism $y:Y\to Z$ of $\scrD$-schemes that extends $y_{\scrK}$.} 

\medskip
A classical result of N\'eron says that each abelian variety over $\scrK$ has a N\'eron model over $\scrD$, cf. [23]. This result has an analogue for the case of torsors of smooth group schemes over $\scrK$ of finite type, cf. [3], Subsection 6.5, Corollary 4. On [3], page 15 it is stated that the importance of the notion of N\'eron models ``seems to be restricted" to ``torsors under group schemes". It was a deep insight of Milne which implicitly pointed out that N\'eron models are important in the study of Shimura varieties, cf. [17], Definitions 2.1, 2.2, 2.5, and 2.9. In this Subsection we bring to a concrete fruition Milne's insight: we will use Corollary 4.3 and [30] to provide large classes of  projective varieties over certain $\scrK$'s which have projective N\'eron models and which often do not admit finite maps into abelian varieties over $\scrK$. For the rest of the paper we will use the notations of Section 1 and Subsections 2.4 and 2.5. 

\medskip\noindent
{\bf 4.4.1. Proposition.} {\it We assume that the principally polarized abelian scheme $(A,\lambda_A)$ over the number field $E$ is such that $A$ has compact factors. We also assume that the reflex field $E_A$ of $(H_A,X_A)$ is unramified at all primes not dividing $N$ and that the $O_{E_A}[{1\over N}]$-scheme $\scrN_N$ of Subsection 2.4 is smooth. Then $\scrN_N$ is the N\'eron model of $\scrN_{NE_A}=\Sh(H_A,X_A)/K_A(N)$ over $O_{E_A}[{1\over N}]$.}

\medskip
\proof
Let $Y$ be a smooth $O_{E_A}[{1\over N}]$-scheme. Let $y_{E_A}:Y_{E_A}\to\scrN_{NE_A}$ be a morphism of $E_A$-schemes. Let $U$ be an open subscheme of $Y$ such that it contains $Y_{E_A}$ and $y_{E_A}$ extends uniquely to a morphism $y_U:U\to\scrN_N$. As the $O_{E_A}[{1\over N}]$-scheme $\scrN_N$ is projective (cf. Corollary 4.3), we can assume that the codimension of $Y\setminus U$ in $Y$ is at least $2$. 

Let $(B_U,\lambda_{B_U}):=y_U^*(\scrB,\lambda_{\scrB})$. The abelian scheme $B_U$ extends to an abelian scheme $B_Y$ over $Y$ (cf. [30], Theorem 1.3) in a unique way (cf. [17], Corollary 2.12). Also $\lambda_{B_U}$ extends uniquely to a principal polarization $\lambda_{B_Y}$ of $B_Y$, cf. [17], Proposition 2.14. Obviously, the level-$N$ symplectic similitude structure of $(B_U,\lambda_{B_U})$ extends uniquely to a level-$N$ symplectic similitude structure of $(B_Y,\lambda_{B_Y})$. Thus the composite of $y_U$ with the finite morphism $\scrN_N\to\scrA_{g,1,N}$ extends uniquely to a morphism $z:Y\to\scrA_{g,1,N}$. As $Y$ is normal and as the morphism $\scrN_N\to\scrA_{g,1,N}$ is finite, $z$ factors uniquely through a morphism $y:Y\to\scrN_N$. Obviously $y$ extends $y_U$ and thus also $y_{E_A}$. From this and the uniqueness of $y$ and $y_U$, we get that $\scrN_N$ satisfies the N\'eron mapping property. Thus $\scrN_N$ is the N\'eron model of $\scrN_{NE_A}=\Sh(H_A,X_A)/K_A(N)$ over $O_{E_A}[{1\over N}]$.\endproof

\medskip\noindent
{\bf 4.4.2. Remark.} If $N$ has many prime divisors, then $K_A(N)$ is a sufficiently small compact open subgroup of $H_A(\dbA_f)$ and thus $\Sh(H_A,X_A)_{\dbC}/K_A(N)$ is a projective, smooth $\dbC$-scheme of general type (see [16], \S2, Subsection 1.2). Thus $\scrN_N$ is not among the N\'eron models studied in [3]. If the Albanese variety of each connected component $\scrC$ of $\Sh(H_A,X_A)_{\dbC}/K_A(N)$ is trivial, then $\Sh(H_A,X_A)/K_A(N)$ is not a finite scheme over an abelian variety over $E_A$. Example: if $H^{\ad}_{A\dbR}\arrowsim\text{SU}(a,b)^{\ad}_{\dbR}\times_{\dbR} \text{SU}(a+b,0)^{\ad}_{\dbR}$, with $a$, $b\in\dbN\setminus\{1,2\}$, then we have $H^{1,0}(\scrC(\dbC),\dbC)=0$ (cf. [25], Theorem 2, 2.8 (i)) and thus the Albanese variety of $\scrC$ is trivial; therefore the connected components of the projective $E_A$-scheme $\Sh(H_A,X_A)/K_A(N)$ are not finite schemes over torsors of smooth groups over $E_A$. This remark was hinted at in [30].

\bigskip\noindent
{\bf 4.5. Example.} We assume that $A$ has compact factors, that $N\in 6\dbN$, and that the Zariski closure $H_{A,N}$ of $H_A$ in $\text{GL}_{L_A[{1\over N}]}$ is a reductive group scheme over $\dbZ[{1\over N}]$. We also assume that if $(H_A^{\ad},X_A^{\ad})$ has a simple factor of $A_n$ type, then either all prime factors of $n+1$ divide $N$ or the degree of the isogeny $H_A^{\sc}\to H_A^{\der}$ divides $N$.${}^1$ $\vfootnote{1} {This condition is not truly needed. It is inserted only to avoid the error which was made in the b) part of [29], Theorem 6.2.2 and which is eliminated in [32].}$ Let $p\in\dbN$ be an arbitrary prime that does not divide $N$; thus $p\Ge 5$. Let $\dbZ_{(p)}$ be the localization of $\dbZ$ with respect to $p$. As ${H_{A,N}}_{\dbZ_{(p)}}$ is a reductive group scheme over $\dbZ_{(p)}$, the field $E_A$ is unramified over $p$ (cf. [18], Corollary 4.7 (a)). Thus the normalization $E_{A(p)}$ of $\dbZ_{(p)}$ in $E_A$ is a finite, \'etale $\dbZ_{(p)}$-algebra. This implies that $E_A$ is unramified at all primes not dividing $N$. In the next two paragraphs we check that $\scrN_N$ is a smooth $O_{E_A}[{1\over N}]$-scheme. 

The $E_A$-scheme $\Sh(H_A,X_A)/K_A(N)$ is smooth and the natural quotient morphism $\Sh(H_A,X_A)\twoheadrightarrow\Sh(H_A,X_A)/K_A(N)$ is a pro-\'etale cover, cf. Subsection 2.4. We define $\scrM:=\text{proj.}\,\text{lim.}_{\tilde N\in N\dbN,g.c.d.(\tilde N,p)=1} {\scrA_{g,1,\tilde N}}_{\dbZ_{(p)}}$. It is well known that we can identify $\scrM_{\dbQ}=\Sh(\text{GSp}(W_A,\psi_A),S_A)/\text{GSp}(L_A,\psi_A)(\dbZ_p)$ and that $\scrM$ is the integral canonical model of the Shimura triple $(\text{GSp}(W_A,\psi_A),S_A,\text{GSp}(L_A,\psi_A)(\dbZ_p))$ as defined in [29], Subsubsections 3.2.6 and 3.2.3 6) (see [17], Theorem 2.10 and [29], Example 3.2.9). Let $\scrN^{(p)}:=\scrN_A^{(p)}$ be the normalization of $\scrM_{E_{A(p)}}$ in (the ring of fractions of) $\Sh(H_A,X_A)/H_{A,N}(\dbZ_p)$; the $E_A$-scheme $\scrN^{(p)}_{E_A}=\Sh(H_A,X_A)/H_{A,N}(\dbZ_p)$ is a pro-\'etale cover of $\scrN_{NE_A}=\Sh(H_A,X_A)/K_A(N)$. 

From proof of [29], Proposition  3.4.1 we get that $\scrN^{(p)}$ is a pro-\'etale cover of $\scrN_{NE_{A(p)}}$ (the previous paragraph implies that conditions (i) and (ii) of loc. cit. hold in the context of $\scrM$ and $\scrN^{(p)}$). As $p\Ge 5$, from [29], Subsubsections 3.4.1, 3.2.12, and 6.4.1 we get that $\scrN^{(p)}$ is the integral canonical model of the Shimura triple $(H_A,X_A,H_{A,N}(\dbZ_p))$. Thus $\scrN^{(p)}$ is a regular, formally smooth $E_{A(p)}$-scheme. This implies that $\scrN_{NE_{A(p)}}$ is a smooth $E_{A(p)}$-scheme. As $p\in\dbN$ was an arbitrary prime that does not divide $N$, we conclude that $\scrN_N$ is a smooth $O_{E_A}[{1\over N}]$-scheme. 

Thus $\scrN_N$ is a N\'eron model of $\scrN_{NE_A}$ over $O_{E_A}[{1\over N}]$, cf. Proposition 4.4.1.

\medskip\noindent
{\bf 4.6. Remarks.} {\bf (a)} Either [14], \S5 or [32] can be used to provide many examples similar to the one of Example 4.5 but with $N$ relatively prime to either $2$ or $3$.

\smallskip
{\bf (b)} We refer to Example 4.5; thus the prime $p$ is at least $5$. As $\scrN_N$ is a projective, smooth $O_{E_A}[{1\over N}]$-scheme (cf. Corollary 4.3 and Example 4.5), $\scrN^{(p)}$ is a pro-\'etale cover of a projective, smooth $E_{A(p)}$-scheme. This validates the erroneous [29], Remark 6.4.1.1 2) for the case of Shimura pairs $(G,X)$ of abelian type that have compact factors. In other words, if the group $G_{\dbQ_p}$ is unramified, then the scheme $\Sh_p(G,X)$ proved to exist in [29], Theorem 6.4.1 is a pro-\'etale cover of a smooth, projective scheme over the normalization  of $\dbZ_{(p)}$ in $E(G,X)$. Based on [29], Subsection 6.8 and Subsubsections 6.8.1 and 6.8.2 a), in the last sentence one can replace ``abelian type" by ``preabelian type". Implicitly, this validates [29], Subsubsection 6.4.11 for all Shimura pairs $(G,X)$ of preabelian type that have compact factors. We recall that a Shimura pair $(G,X)$ is said to be of preabelian type, if $(G^{\ad},X^{\ad})$ is isomorphic to $(H_A^{\ad},X_A^{\ad})$ for some abelian variety $A$ over a number field $E$ (see [17], [29], etc.). If moreover one can assume that we have central isogenies $H_A^{\der}\to G^{\der}\to G^{\ad}\arrowsim H_A^{\ad}$, then $(G,X)$ is said to be of abelian type.

\bigskip\smallskip
\references{37}
{\nspace{

\Ref[1]
W. Baily and A. Borel,
\sl Compactification of arithmetic quotients of bounded
symmetric domains,
\rm Ann. of Math. (2) {\bf 84} (1966), no. 3, pp. 442--528.

\Ref[2]
A. Borel and Harish-Chandra,
\sl Arithmetic subgroups of algebraic groups,
\rm Ann. of Math. (2) {\bf 75} (1962), no. 3, pp. 485--535.

\Ref[3]
S. Bosch, W. L\"utkebohmert, and M. Raynaud,
\sl N\'eron models,
\rm Ergebnisse der Mathematik und ihrer Grenzgebiete (3), Vol. {\bf 21}, Springer-Verlag, Berlin, 1990.

\Ref[4] N. Bourbaki, 
\sl Lie groups and Lie algebras, Chapters 4--6,
\rm Elements of Mathematics (Berlin), Springer-Verlag, Berlin, 2002.

\Ref[5] P. Berthelot and A. Ogus, 
\sl F-crystals and de Rham cohomology. I, 
\rm Invent. Math. {\bf 72} (1983), no. 2, pp. 159--199.

\Ref[6]
P. Deligne,
\sl Travaux de Shimura,
\rm S\'eminaire  Bourbaki, 23\`eme ann\'ee (1970/71), Exp. No. 389,  pp. 123--165, Lecture Notes in Math., Vol. {\bf 244}, Springer-Verlag, 1971.

\Ref[7]
P. Deligne,
\sl Vari\'et\'es de Shimura: interpr\'etation modulaire, et
techniques de construction de mod\`eles canoniques,
\rm Automorphic forms, representations and $L$-functions (Oregon State Univ., Corvallis, OR, 1977), Part 2,  pp. 247--289, Proc. Sympos. Pure Math., Vol. {\bf 33}, Amer. Math. Soc., Providence, RI, 1979.

\Ref[8]
P. Deligne,
\sl Hodge cycles on abelian varieties,
\rm Hodge cycles, motives, and Shimura varieties, Lecture Notes in Math., Vol. {\bf 900}, pp. 9--100, Springer-Verlag, 1982.

\Ref[9] 
M. Demazure, A. Grothendieck, et al., 
\sl Sch\'emas en groupes, Vol. II, 
\rm Lecture Notes in Math., Vol. {\bf 152}, Springer-Verlag, 1970.

\Ref[10]
G. Faltings and C.-L. Chai,
\sl Degeneration of abelian varieties,
\rm  Ergebnisse der Mathematik und ihrer Grenzgebiete (3), Vol. {\bf 22}, Springer-Verlag, Berlin, 1990.

\Ref[11] 
J.-M. Fontaine, 
\sl Repr\'esentations p-adiques semi-stables,
\rm P\'eriodes $p$-adiques (Bures-sur-Yvette, 1988), J. Ast\'erisque {\bf 223}, pp. 113--184, Soc. Math. de France, Paris, 1994.

\Ref[12] 
A. Grothendieck, 
\sl \'El\'ements de g\'eom\'etrie alg\'ebrique. IV. \'Etude locale des sch\'emas et des morphismes de sch\'ema (Quatri\`eme Partie), 
\rm Inst. Hautes \'Etudes Sci. Publ. Math., Vol. {\bf 32}, 1967.

\Ref[13]
R. Hartshorne,
\sl Algebraic geometry,
\rm Grad. Texts in Math. {\bf 52}, Springer-Verlag, New York--Heidelberg, 1977.

\Ref[14]
R. E. Kottwitz,
\sl Points on some Shimura Varieties over finite fields,
\rm J. of Amer. Math. Soc. {\bf 5} (1992), no. 2, pp. 373--444.

\Ref[15]
H. Matsumura,
\sl Commutative algebra. Second edition,
\rm The Benjamin/Cummings Publishing Co., Inc., 1980.

\Ref[16]
J. S. Milne,
\sl Canonical models of (mixed) Shimura varieties and automorphic vector bundles, 
\rm  Automorphic Forms, Shimura varieties and L-functions, Vol. I (Ann Arbor, MI, 1988), pp. 283--414, Perspectives in Math., Vol. {\bf 10}, Acad. Press, Boston, MA, 1990.

\Ref[17]
J. S. Milne, 
\sl The points on a Shimura variety modulo a prime of good
reduction,
\rm The Zeta function of Picard modular surfaces, pp. 153--255, Univ. Montr\'eal, Montreal, Quebec, 1992.

\Ref[18]
J. S. Milne,
\sl Shimura varieties and motives,
\rm  Motives (Seattle, WA, 1991), Part 2, pp. 447--523, Proc. Symp. Pure Math., Vol. {\bf 55}, Amer. Math. Soc., Providence, RI, 1994.

\Ref[19]
J. S. Milne,
\sl Descent for Shimura varieties,
\rm Mich. Math. J. {\bf 46} (1999), no. 1, pp. 203--208.

\Ref[20]
Y. Morita,
\sl On potential good reduction of abelian varieties,
\rm J. Fac. Sci. Univ. Tokyo Sect. I A Math. {\bf 22} (1975), no. 3, pp. 437--447.

\Ref[21]
D. Mumford, J. Fogarty, and F. Kirwan,
\sl Geometric invariant theory. Third enlarged edition, 
\rm Ergebnisse der Mathematik und ihrer Grenzgebiete (2), Vol. {\bf 34}, Springer-Verlag, Berlin, 1994. 

\Ref[22]
D. Mumford,
\sl Abelian varieties,
\rm Tata Inst. of Fund. Research Studies in Math., No. {\bf 5}, Published for the Tata Institute of Fundamental Research, Bombay; Oxford Univ. Press, London, 1970 (reprinted 1988).

\Ref[23] 
A. N\'eron, 
\sl Mod\`eles minimaux des vari\'et\'es ab\'eliennes, 
\rm Inst. Hautes \'Etudes Sci. Publ. Math., Vol. {\bf 21}, 1964.

\Ref[24]
R. Noot,
\sl Abelian varieties with $l$-adic Galois representations of Mumford's type,
\rm J. reine angew. Math. {\bf 519} (2000), pp. 155--169.

\Ref[25]
R. Parthasarathy, 
\sl Holomorphic forms in $\Gamma\backslash G/K$ and Chern classes,
\rm Topology {\bf 21}  (1982), no. 2, pp. 157--178.

\Ref[26]
F. Paugam,
\sl Galois representations, Mumford--Tate groups and good reduction of abelian varieties,
\rm Math. Ann. {\bf 329} (2004), no. 1, pp. 119--160. Erratum: Math. Ann. {\bf 332} (2004), no. 4, p. 937. 

\Ref[27]
T. A. Springer,
\sl Reductive groups,
\rm Automorphic forms, representations and $L$-functions (Oregon State Univ., Corvallis, Ore., 1977), Part 1, pp. 3--27, Proc. Sympos. Pure Math., Vol. {\bf 33}, Amer. Math. Soc., Providence, RI, 1979.

\Ref[28] 
J. Tits, 
\sl Classification of algebraic semisimple groups, 
\rm Algebraic Groups and Discontinuous Subgroups (Boulder, CO, 1965),  pp. 33--62, Proc. Sympos. Pure Math., Vol. {\bf 9}, Amer. Math. Soc., Providence, RI, 1966.

\Ref[29]
A. Vasiu,
\sl Integral canonical models of Shimura varieties of preabelian type,
\rm Asian J. Math. {\bf 3} (1999), no. 2, pp. 401--518.

\Ref[30]
A. Vasiu,
\sl A purity theorem for abelian schemes,
\rm Mich. Math. J. {\bf 54} (2004), no. 1, pp. 71--81. 

\Ref[31]
A. Vasiu,
\sl The Mumford--Tate Conjecture and Shimura Varieties, Part I,
\rm math.NT/0212066.

\Ref[32]
A. Vasiu,
\sl Integral canonical models of unitary Shimura varieties, 
\rm math.NT/0505507.

}}

\bigskip
\hbox{Adrian Vasiu,\;\;\;Email: adrian\@math.arizona.edu}
\hbox{Address: University of Arizona, Department of Mathematics, 617 N. Santa Rita Ave.,}
\hbox{P.O. Box 210089, Tucson, AZ 85721-0089, U.S.A.}

\enddocument